\documentclass[preprint,12pt]{elsarticle}

\usepackage{amsmath,amssymb,amsthm,bm}

\usepackage{booktabs}
\usepackage[ruled,linesnumbered]{algorithm2e}
\usepackage{hyperref}
\usepackage{float}

\usepackage{graphicx}
\usepackage{caption}
\usepackage{hyperref}

\usepackage{amssymb}
\usepackage{color}
\usepackage{amsmath}
\newcommand{\stor}{\operatorname{stor}}
\usepackage{bm}
\usepackage{subfigure}
\usepackage{float}
\usepackage{hyperref}
\usepackage[capitalize]{cleveref}
\usepackage{booktabs}
\usepackage{mathtools}
\makeatletter
\def\ps@pprintTitle{%
  \let\@oddhead\@empty
  \let\@evenhead\@empty
  \let\@oddfoot\@empty
  \let\@evenfoot\@empty
}
\makeatother

\newcommand{\rank}{\operatorname{rank}}

\newtheorem{theorem}{Theorem}[section]
\journal{Journal of Computational Physics}
\begin{document}

\begin{frontmatter}

\title{An Isogeometric Collocation Solver Based on the Quantized Tensor Train Representation}

\author[label1]{Jun Song}
 \ead{junsong@whu.edu.cn}
 
\author[label1,label2]{Hua Xiang\corref{cor1}} 
\ead{hxiang@whu.edu.cn}
\cortext[cor1]{Corresponding author.}

\affiliation[label1]{organization={School of Mathematics and Statistics, Wuhan University},
            city={Wuhan},
            postcode={430072}, 
            country={China}}

\affiliation[label2]{organization={ Hubei Key Laboratory of Computational Science, Wuhan University},
            city={Wuhan},
            postcode={430072}, 
            country={China}}

\begin{abstract}
  We introduce an isogeometric collocation framework based on the quantized
tensor train (QTT) representation to solve partial differential equations. 
The proposed method constructs isogeometric collocation operators directly in
the QTT format using only a few sampled entries, enabling efficient compression and
computation.
The  vectors of the
geometry-derived coefficients at the collocation points are represented
in the QTT format using the TT-cross approximation. 
We conduct numerical experiments using the proposed QTT-based isogeometric
collocation framework for the convection-diffusion-reaction problems in different
spatial dimensions and for the Poisson problems on three-dimensional nontrivial geometries,
showing significant reductions in memory usage and computational cost while
preserving solution accuracy.
\end{abstract}


\end{frontmatter}

\section{Intoduction}

High-accuracy numerical simulation of high-dimensional partial differential equations  (PDEs) is of crucial importance in engineering and scientific applications. In recent years, tensor networks (TNs) \cite{bachmayr2016} have emerged as a promising strategy to eliminate the curse of dimensionality. Representative tensor factorizations include the Canonical Polyadic (CP) decomposition, the Hierarchical Tucker (HT) decomposition \cite{kolda2009}, and the tensor train (TT) decomposition \cite{oseledets2011}.
The TT format is proposed as a computationally efficient alternative to the CP and Tucker decompositions, representing a \(d\)-dimensional tensor as a chain of low-rank cores \cite{oseledets2011}. Subsequent advances include the TT-cross approximation, which builds a TT representation from a limited number of tensor entries \cite{oseledets2010} and the quantized TT (QTT) format, which attains logarithmic compression by folding long vectors into higher-order tensors \cite{kazeev2012}. These developments have enabled efficient solutions of high-dimensional PDEs \cite{dolgov2012,khoromskij2011}.

Recent applications of tensor networks to PDEs include efficient and accurate solutions of the Poisson equation \cite{khoromskij2011qtt}, the time-independent Schr\"odinger equation \cite{gidi2025}, Maxwell's equations \cite{manzini2023}, the Fokker--Planck equation \cite{dolgov2012}, the time-dependent convection--diffusion--reaction equation \cite{adak2024}, and incompressible fluid dynamics \cite{peddinti2024}. TN-based techniques have successfully exploited the low-rank structure of differential operators in the finite difference  method \cite{manzini2023}, the spectral method \cite{garcia2021}, the pseudospectral method \cite{gidi2025}, and the spectral collocation method \cite{adak2024}. Most of these applications have been largely confined to problems discretized on regular Cartesian domains. For non-Cartesian domains, recent works have identified  isogeometric analysis (IGA) \cite{cottrell2009,hughes2005}  as a natural candidate for integration with tensor network solvers.

Isogeometric analysis uses tensor-product B-spline bases for both the geometry description and the discretization of partial differential equations. IGA is typically implemented within two frameworks: IGA Galerkin \cite{hughes2005} and IGA collocation \cite{auricchio2010}. When Hughes et al. proposed IGA, the computation was implemented by the Galerkin method \cite{zienkiewicz2013} through the weak form, which involved the expensive multidimensional quadrature and thus significantly increased the complexity and computation time of matrix formation. The isogeometric collocation method \cite{auricchio2010}, proposed by Auricchio et al., selected collocation points in the computational domain and directly dealt with the strong form of the PDE. The convergence and consistency of some theoretical results for IGA collocation have been examined and proven by Lin et al.~\cite{lin2022,lin2020}. IGA Galerkin provides the optimal convergence, while IGA collocation performs better in terms of the ratio of accuracy to computational time \cite{math11020469}.

 So far, low-rank and tensor-structured numerical methods for isogeometric analysis (IGA) have been developed almost exclusively for Galerkin discretizations.  Mantzaflaris et al.~\cite{mantzaflaris2017} introduced an IGA Galerkin framework based on Kronecker-structured assembly.
Markeeva et al.~\cite{markeeva2021}
constructed two-dimensional stiffness matrices using the QTT format on uniform parametric grids.
Ion et al.~\cite{ion2022} 
developed a TT-based  IGA Galerkin  framework  for PDE approximation on parameter dependent geometries
by constructing a low-rank  substitute of the geometry mapping.
Montardini et
al.~\cite{montardini2023} developed an IGA Galerkin solver based on Tucker tensor for
three-dimensional Poisson with truncated preconditioned conjugate gradients. Tran et al.~\cite{tran2026} assembled  the geometry-derived coefficient fields and
discrete operators for large-scale 3D Poisson problems in the TT format via adaptive TT-cross sampling.
Riemer and Stoll~\cite{riemerstoll2025} developed a low-rank adaptive isogeometric framework for Galerkin discretization based on truncated hierarchical B-splines. Overall, these studies establish low-rank
assembly and solution techniques in the discretization of IGA Galerkin.
The tensor-based solution method for IGA collocation  remains much less developed.

In this work, we propose an isogeometric collocation framework based on the QTT representation. We first discuss an efficient and low-cost QTT construction for one-dimensional collocation
matrices. For regular Cartesian domains, we point out that the QTT representation of the
\(d\)-dimensional discrete operator can be constructed using the assembly method
in \cite{kazeev2012}. For Poisson problems on
3D nontrivial geometries, we
reformulate the problem equivalently on the parameter domain, so that all
geometric effects are contained in a set of variable coefficients. The
geometry-derived coefficient vectors at the collocation points are represented
in the QTT format using TT-cross and then combined with the QTT representations of
the one-dimensional collocation matrices to assemble the QTT representation of
the full operator.

The paper is organized as follows. Section~2 reviews the fundamentals of the
tensor train and isogeometric collocation. Section~3 presents an efficient
QTT construction for one-dimensional collocation matrices and, on regular Cartesian domains, derives a QTT representation of the
\(d\)-dimensional IGA collocation operator. Section~4 develops a QTT-based isogeometric
collocation solver for 3D Poisson problems under 3D nontrivial geometries. Section~5 reports numerical results for verifying the methods presented in
Sections~3 and~4. Section~6 concludes the paper.

$Notations.$ All vectors, matrices and tensors are represented by bold letters. The conjugate transpose is denoted by
\((\cdot)^{\mathrm T}\). The norm \(\|\cdot\|\) denotes the Euclidean
norm of a vector or the spectral norm of a matrix, while \(\|\cdot\|_F\)
denotes the Frobenius norm. The symbol \(|\cdot|\) denotes the absolute value of a scalar or the number of elements in a set or list. The set \(\{1,2,\ldots,j\}\) is denoted by
\([j]\).

\section{Preliminary}

\subsection{Tensor Train Representation }

This subsection introduces the tensor train  format and its basic
algebraic operations. Let $\boldsymbol{\mathbf X}$ be a $d$-dimensional
$n_1\times\cdots\times n_d$-vector. If there exist  TT cores
$\boldsymbol{X}_k\in\mathbb{R}^{r_{k-1}\times n_k\times r_k}$,
$k=1,\ldots,d$, with $r_0=r_d=1$, depending on the
corresponding mode index $i_k$ and one or two rank indices, such that
\begin{equation}
\boldsymbol{\mathbf X}(i_1\ldots i_d)
=\sum_{\alpha_1=1}^{r_1}\cdots
\sum_{\alpha_{d-1}=1}^{r_{d-1}}
\boldsymbol{X}_1(1,i_1,\alpha_1)\cdot \boldsymbol{X}_2(\alpha_1,i_2,\alpha_2)
\cdots \boldsymbol{X}_d(\alpha_{d-1},i_d,1),
\label{eq:tt_vector}
\end{equation}
then $\boldsymbol{\mathbf X}$ is said to be represented in the TT-vector format, where $i_k=1,\ldots,n_k$ is called the mode induce  and the summation limit
$r_k$ of the rank induce $\alpha_k$ is called  the TT rank of the decomposition. Let
$n=\max_k n_k$ and $r=\max_k r_k$. The storage cost is then
$\mathcal{O}(dnr^2)$. In quantum physics, a TT vector is commonly
called a matrix product state.

The TT format can also be applied to a $d$-dimensional matrix
$\boldsymbol{\mathbf W}$ of size
$(m_1\times\cdots\times m_d)\times(n_1\times\cdots\times n_d)$. If there exist
TT cores
$\boldsymbol{W}_k\in\mathbb{R}^{R_{k-1}\times m_k\times n_k\times R_k}$,
with $R_0=R_d=1$, such that
\begin{align}
    \boldsymbol{\mathbf W}
    \bigl((i_1\ldots i_d),(j_1\ldots j_d)\bigr)
    ={}&
    \sum_{\gamma_1=1}^{R_1}\cdots
    \sum_{\gamma_{d-1}=1}^{R_{d-1}}
    \boldsymbol{W}_1(1,i_1,j_1,\gamma_1)
    \cdot \boldsymbol{W}_2(\gamma_1,i_2,j_2,\gamma_2)
    \notag\\
    &\cdots
    \boldsymbol{W}_d(\gamma_{d-1},i_d,j_d,1),
    \label{eq:tt_matrix}
\end{align}
then $\boldsymbol{\mathbf W}$ is represented in the TT-matrix format, where
$i_k=1,\ldots,m_k$, $j_k=1,\ldots,n_k$ are the mode induces, and  the summation limit $R_k$ of the rank induce $\gamma_k$  is the corresponding TT
rank.  Let
$m=\max_k m_k$, $n=\max_k n_k$, and $R=\max_k R_k$. The storage cost is
$\mathcal{O}(dmnR^2)$. In quantum physics, the corresponding representation
is commonly called a matrix product operator.

\subsubsection{Basic Algebraic Operations in TT}

We first introduce the core matrix and the rank core product~\cite{kazeev2012}.  Consider a TT core 
$\boldsymbol{W}\in\mathbb{R}^{p\times a\times b\times q}$ with two rank indices $\alpha=1,\ldots,p,\,\beta=1,\ldots,q$, and two mode indices $i=1,\ldots,a,\,j=1,\ldots,b$. Fixing its left
and right rank indices $\alpha$ and $\beta$,
respectively, we define the matrix
$\boldsymbol{A}_{\alpha\beta}
=\boldsymbol{W}(\alpha,:,:,\beta)\in\mathbb{R}^{a\times b}$.
Then
the core $\boldsymbol{W}$ can be considered as the matrix
\[
\begin{bmatrix}
\boldsymbol{A}_{11}&\cdots&\boldsymbol{A}_{1q}\\
\vdots&&\vdots\\
\boldsymbol{A}_{p1}&\cdots&\boldsymbol{A}_{pq}
\end{bmatrix},
\]
which is referred to as the core matrix of $\boldsymbol{W}$.

Let the adjacent cores $\boldsymbol{W}_1$ and $\boldsymbol{W}_2$ have ranks
$R_0\times R_1$ and $R_1\times R_2$, respectively, and consist of blocks
$\boldsymbol{A}^{(1)}_{\alpha_0\alpha_1}\in\mathbb{R}^{m_1\times n_1}$ and
$\boldsymbol{A}^{(2)}_{\alpha_1\alpha_2}\in\mathbb{R}^{m_2\times n_2}$, where
$\alpha_k=1,\ldots,R_k$ for $k=0,1,2$. Let us define a rank core product $\boldsymbol{W}_1\bowtie \boldsymbol{W}_2$ of $\boldsymbol W_1$
and $\boldsymbol W_2$ as a core of rank $R_0\times R_2$, consisting of blocks
\begin{equation}
\boldsymbol{A}_{\alpha_0\alpha_2}
=
\sum_{\alpha_1=1}^{R_1}
\boldsymbol{A}^{(1)}_{\alpha_0\alpha_1}
\otimes \boldsymbol{A}^{(2)}_{\alpha_1\alpha_2}.
\label{eq:rank_core_product}
\end{equation}
of mode size $m_1m_2 \times n_1n_2$.
The TT vector and matrix in \eqref{eq:tt_vector} and \eqref{eq:tt_matrix} can
therefore be recast as
$\boldsymbol{\mathbf X}=\boldsymbol{X}_1\bowtie\cdots\bowtie\boldsymbol{X}_d$
and
$\boldsymbol{\mathbf W}=\boldsymbol{W}_1\bowtie\cdots\bowtie\boldsymbol{W}_d$,
respectively.

We next consider the sum of two tensors in the TT format. Let
$\boldsymbol{\mathbf U}
=\boldsymbol{U}_1\bowtie\cdots\bowtie\boldsymbol{U}_d$ and
$\boldsymbol{\mathbf H}
=\boldsymbol{H}_1\bowtie\cdots\bowtie\boldsymbol{H}_d$
have the same mode size, with TT ranks $r_k$ and $s_k$, respectively. Their
sum can be represented as
\begin{equation}
\boldsymbol{\mathbf U}+\boldsymbol{\mathbf H}
=
\begin{bmatrix}\boldsymbol{U}_1&\boldsymbol{H}_1\end{bmatrix}
\bowtie
\begin{bmatrix}\boldsymbol{U}_2& \boldsymbol 0\\\boldsymbol 0&\boldsymbol{H}_2\end{bmatrix}
\bowtie\cdots\bowtie
\begin{bmatrix}\boldsymbol{U}_{d-1}&\boldsymbol 0\\ \boldsymbol 0&\boldsymbol{H}_{d-1}\end{bmatrix}
\bowtie
\begin{bmatrix}\boldsymbol{U}_d\\\boldsymbol{H}_d\end{bmatrix}.
\label{eq:tt_sum}
\end{equation}
Here addition is meant
elementwise and the square brackets denote core matrices.  The resulting TT ranks
are at most $r_k+s_k$, $k=1,\ldots,d-1$.

The tensor product of two tensors in the TT format is
\begin{equation}
\boldsymbol{\mathbf U}\otimes\boldsymbol{\mathbf H}
=\boldsymbol{U}_1\bowtie\cdots\bowtie\boldsymbol{U}_d
\bowtie\boldsymbol{H}_1\bowtie\cdots\bowtie\boldsymbol{H}_d.
\label{eq:tt_Kronecker}
\end{equation}
Since the connection
rank is equal to $1$, no additional rank growth is introduced.

The matrix-vector product can be represented
as
\[\boldsymbol{\mathbf W}\boldsymbol{\mathbf X}
=(\boldsymbol{W}_1\star\boldsymbol{X}_1)
\bowtie\cdots\bowtie
(\boldsymbol{W}_d\star\boldsymbol{X}_d)\],
where $\star$ denotes the contraction of the corresponding cores over their
common mode index. This operation is defined by
\begin{equation}
\begin{aligned}
&\bigl(\boldsymbol{W}_k\star\boldsymbol{X}_k\bigr)
\bigl((\gamma_{k-1},\alpha_{k-1}),i_k,(\gamma_k,\alpha_k)\bigr)
=
\sum_{j_k=1}^{n_k}
\boldsymbol{W}_k(\gamma_{k-1},i_k,j_k,\gamma_k)
\boldsymbol{X}_k(\alpha_{k-1},j_k,\alpha_k).
\end{aligned}
\label{eq:tt_local_matvec}
\end{equation}
Therefore, the resulting core has left and right rank dimensions
$R_{k-1}r_{k-1}$ and $R_kr_k$, respectively.

The above algebraic operations may increase the TT ranks. In practical
computations, TT rounding is performed with a prescribed tolerance
 in order to control the storage cost and computational
complexity.

\subsubsection{Quantized Tensor Train}

 Let the  $d$-dimensional
$n_1\times\cdots\times n_d$ vector
${\boldsymbol{\mathbf X}}$ satisfy
$n_k=2^{q_k}$ with $q_k\in\mathbb{N}$. Expand each mode index
$i_k=1,\ldots,n_k$ as
\[
i_k=1+\sum_{\ell=1}^{q_k}(i_{k,\ell}-1)2^{q_k-\ell},
\]
where $i_{k,\ell}\in\{1,2\}$. The original tensor can then be reshaped into a
binary tensor of order $M=\sum_{k=1}^d q_k=\log_2 N$, where
$N=\prod_{k=1}^d n_k$. Denote all indices of virtual modes $\{i_{1,1},\ldots,i_{1,q_1},i_{2,1},\ldots,i_{2,q_2},\ldots,i_{d,1},\ldots,i_{d,q_d}\}$ consecutively by
$\widehat{i}_1,\ldots,\widehat{i}_M$. The QTT-vector representation of $\boldsymbol{\mathbf X}$ is
\[
{\boldsymbol{\mathbf X}}
=
\widehat{\boldsymbol{X}}_1\bowtie\cdots
\bowtie\widehat{\boldsymbol{X}}_M,
\]
where
$\widehat{\boldsymbol{X}}_\ell
\in\mathbb{R}^{\widehat{r}_{\ell-1}\times 2\times\widehat{r}_\ell}$
and $\widehat{r}_0=\widehat{r}_M=1$. The ranks
$\widehat{r}_1,\ldots,\widehat{r}_{M-1}$ are called QTT ranks of the
representation. Let $\widehat{r}=\max_\ell\widehat{r}_\ell$, which is hereafter denoted as $\operatorname{rank}_{\mathrm{QTT}}({\boldsymbol{\mathbf X}})$.
With $n=\max_k n_k$, the storage cost of a QTT vector is
$\mathcal{O}(d\widehat{r}^2\log_2 n)$.

For the $d$-dimensional matrix
$\boldsymbol{\mathbf W}$ of size
$(m_1\times\cdots\times m_d)\times(n_1\times\cdots\times n_d)$, consider the square case $m_k=n_k=2^{q_k}$. Expand the row
and column mode indices $i_k$ and $j_k$ in the same way, and
denote their indices of virtual modes consecutively by $\widehat{i}_\ell$ and
$\widehat{j}_\ell$.  The
QTT-matrix representation of ${\boldsymbol{\mathbf W}}$ is
\[
{\boldsymbol{\mathbf W}}
=
\widehat{\boldsymbol{W}}_1\bowtie\cdots
\bowtie\widehat{\boldsymbol{W}}_M,
\]
where
$\widehat{\boldsymbol{W}}_\ell
\in\mathbb{R}^{\widehat{R}_{\ell-1}\times 2\times 2\times\widehat{R}_\ell}$
and $\widehat{R}_0=\widehat{R}_M=1$. The ranks
$\widehat{R}_1,\ldots,\widehat{R}_{M-1}$ are QTT-matrix ranks. Let
$\widehat{R}=\max_\ell\widehat{R}_\ell$, which is hereafter denoted as
$\rank_{\mathrm{QTT}}({\boldsymbol{\mathbf W}})$ . The storage cost of the QTT-matrix is
$\mathcal{O}(d\widehat{R}^2\log_2 n)$.

Consequently, when the QTT ranks are bounded independently of the number of
virtual dimensions, the QTT format reduces the polynomial dependence of the
storage cost on the mode size $n$ to a logarithmic dependence. This
substantially reduces the storage cost for large-scale discretized problems.

\subsection{Isogeometric Collocation Method}

Consider the boundary value problem on the physical domain \(\Omega\subset\mathbb R^d\),
\begin{equation}
\begin{cases}
\mathcal L u(\boldsymbol x)=f(\boldsymbol x), & \boldsymbol x\in\Omega,\\[2pt]
\mathcal B u(\boldsymbol x)=g(\boldsymbol x), & \boldsymbol x\in\partial\Omega,
\end{cases}
\label{eq:iga_bvp}
\end{equation}
where $\Omega$ is an open domain, $\partial\Omega$ is its boundary, $u$ is the function to be solved, $\mathcal L$ is a differential operator, $\mathcal B$ is a boundary operator, and $f$ and $g$ are known functions.

Assume
\(
\{\widehat\varphi_{\mathbf i}(\boldsymbol\xi):
\boldsymbol\xi\in\widehat{\Omega},\ \mathbf i\in\mathcal I\}\) is a set of basis functions defined on the parameter domain \(\widehat{\Omega}:=[0,1]^d\), where the index set is
\(
\mathcal I
=
\{1,\ldots,n_1\}\times\cdots\times\{1,\ldots,n_d\},
\mathbf i=(i_1,\ldots,i_d).
\)
Suppose that the physical domain \(\Omega\) can be described by a geometry mapping function \(\boldsymbol{G}:\widehat{\Omega}\to\overline{\Omega}\), \(\boldsymbol{x} = \boldsymbol{G}(\boldsymbol{\xi})\). Here,  \(\overline{\Omega}=\Omega\cup\partial\Omega\). Under the standard isogeometric analysis format, this mapping can be expressed as
\begin{equation}
\boldsymbol G(\boldsymbol\xi)
=
\sum_{\mathbf i\in\mathcal I}
\widehat\varphi_{\mathbf i}(\boldsymbol\xi)\,
\boldsymbol P_{\mathbf i},
\qquad
\boldsymbol P_{\mathbf i}\in\mathbb R^d,
\label{eq:iga_geometry}
\end{equation}
where \(\boldsymbol P_{\mathbf i}\) represents the control point. Assume that \(\boldsymbol G\) is invertible on \(\widehat{\Omega}\) and satisfies \(\boldsymbol G(\partial\widehat{\Omega})=\partial\Omega\). In the most commonly used format of isogeometric analysis, \(\widehat\varphi_{\mathbf i}\) is taken as the B-spline or NURBS basis function. Nevertheless, other types of functions can also be used.

Let \(u_h(\boldsymbol x)\) denote the numerical solution in the physical domain.
From the geometry mapping \(\boldsymbol x=\boldsymbol G(\boldsymbol\xi)\), we have \(u_h(\boldsymbol x)=u_h(\boldsymbol G(\boldsymbol{\xi}))\). Define
\(
\widehat u_h(\boldsymbol\xi)
=
u_h\bigl(\boldsymbol G(\boldsymbol\xi)\bigr).
\) We approximate the parametric-domain numerical solution using B-spline or NURBS basis functions, 
\begin{equation}
\widehat u_h(\boldsymbol\xi)
=
\sum_{\mathbf i\in\mathcal I}
c_{\mathbf i}\widehat\varphi_{\mathbf i}(\boldsymbol\xi),
\label{eq:iga_parametric_solution}
\end{equation}
where \(c_{\mathbf i}\) are the unknown coefficients to be determined.

A set of collocation points is introduced in the parameter domain \(\widehat{\Omega}\),
\begin{equation}
  \{\bar{\boldsymbol\xi}_{\mathbf i}\in\widehat{\Omega}:\ \mathbf i\in\mathcal I\}.  \label{eq:collocation points}
\end{equation}
The index set \(\mathcal I\) is divided into the interior collocation index set \(\mathcal I^{(\mathrm{int})}\) and the boundary collocation index set \(\mathcal I^{(\mathrm{bd})}\), depending on whether the collocation points lie inside the parameter domain.
The corresponding collocation points in the physical domain are obtained via the geometry mapping,
\(
\bar{\boldsymbol x}_{\mathbf i}
=
\boldsymbol G(\bar{\boldsymbol\xi}_{\mathbf i}),
\, \mathbf i\in\mathcal I.
\)

The isogeometric collocation method requires the numerical solution to directly satisfy the governing equations and boundary conditions at all collocation points.
\begin{equation}
\begin{cases}
\mathcal L u_h(\bar{\boldsymbol x}_{\mathbf i})
=
f(\bar{\boldsymbol x}_{\mathbf i}),
& \mathbf i\in\mathcal I^{(\mathrm{int})},\\[3pt]
\mathcal B u_h(\bar{\boldsymbol x}_{\mathbf i})
=
g(\bar{\boldsymbol x}_{\mathbf i}),
& \mathbf i\in\mathcal I^{(\mathrm{bd})}.
\end{cases}
\label{eq:iga_collocation}
\end{equation}
Since the differential operator is defined with respect to the physical coordinates \(\boldsymbol{x}\), whereas the basis functions are defined in the parametric coordinates \(\boldsymbol{\xi}\), the Jacobian matrix of the geometry mapping
\(
\boldsymbol{J}
=
{\partial\boldsymbol{x}}/{\partial\boldsymbol{\xi}}
\)
must be introduced to transform derivatives via the chain rule. The first-order derivatives satisfy
\(
\nabla_{\boldsymbol{x}} u_h
=
\boldsymbol{J}^{-\mathsf{T}}
\nabla_{\boldsymbol{\xi}}\widehat{u}_h.
\)
For differential operators involving second derivatives, repeated application of the chain rule yields the corresponding second derivatives in the physical domain.

We adopt the setting of~\cite{donatelli2015robust}: the basis functions in the discrete solution \eqref{eq:iga_parametric_solution} are taken to be tensor-product B-splines over open uniform knot sequences, and the collocation points are taken to be the Greville abscissae. The geometry mapping \(\boldsymbol{G}\) may be any smooth invertible mapping, not necessarily expressed in terms of \eqref{eq:iga_geometry}. Under this setting, \cite{donatelli2015robust} and \cite{donatelli2015spectral} analyzed the asymptotic eigenvalue distribution of IGA collocation matrices and accordingly developed preconditioners and iterative solvers.

We now give the explicit construction of the basis functions and collocation points. In the \(k\)-th parametric direction, the number of basis functions is \(n_k\) and the degree is \(p_k\). The associated open uniform knot vector \(\Xi^{(k)}=
\{\xi_1^{(k)},\ldots,\xi_{n_k+p_k+1}^{(k)}\}\) is defined by
\begin{align*}
   &\xi_1^{(k)}=\cdots=\xi_{p_k+1}^{(k)}=0,\,\,
\xi_{p_k+1+j}^{(k)}=j h_k,\quad j=1,\ldots,n_k-p_k-1,\,\, \\&
\xi_{n_k+1}^{(k)}=\cdots=\xi_{n_k+p_k+1}^{(k)}=1
\end{align*}
where \(h_k=1/(n_k-p_k)\) denotes the distance between interior knots in the \(k\)-th direction. 
The multiplicity of this nodes at both endpoints are \(p_k + 1\), and the internal nodes are uniformly distributed.

According to the Cox–De Boor recursion formula, the B-spline basis function in the \(k\)-th direction is recursively generated. The zeroth-order basis function is
\[
N^{(k)}_{i,0}(\xi)
=
\begin{cases}
1, & \xi_i^{(k)}\le\xi<\xi_{i+1}^{(k)},\\[2pt]
0, & \text{otherwise}.
\end{cases}
\]
For \(p_k \ge 1\),
\begin{equation}
N^{(k)}_{i,p_k}(\xi)
={}
\frac{\xi-\xi_i^{(k)}}{\xi_{i+p_k}^{(k)}-\xi_i^{(k)}}
N^{(k)}_{i,p_k-1}(\xi)+
\frac{\xi_{i+p_k+1}^{(k)}-\xi}{\xi_{i+p_k+1}^{(k)}-\xi_{i+1}^{(k)}}
N^{(k)}_{i+1,p_k-1}(\xi).
\label{eq:iga_cox_de_boor}
\end{equation}
 The support set of the basis function \(N^{(k)}_{i,p_k}(i=1,\ldots,n_k)\) is
\(
\operatorname{supp}\,\bigl(N^{(k)}_{i,p_k}\bigr)
=
[\xi_i^{(k)},\xi_{i+p_k+1}^{(k)}].
\) 
The first and last knots having multiplicity \(p+1\)  guarantees the endpoint interpolation,
\(
N_{i,p}(0)=\delta_{i0},\,
N_{i,p}(1)=\delta_{in_k},
\)where $\delta_{ij}$ denotes the Kronecker delta.
The Greville abscissa of the \(k\)-th direction is defined as
\[
\bar{\xi}^{(k)}_i
=
\frac{1}{p_k}
\sum_{m=1}^{p_k}\xi_{i+m}^{(k)},
\qquad
i=1,\ldots,n_k.
\]
From the open uniform node vector, we know that \(\bar{\xi}^{(k)}_1 = 0\), \(\bar{\xi}^{(k)}_{n_k} = 1\), and \(\{\bar{\xi}^{(k)}_i\}\) is strictly increasing with respect to \(i\). Therefore, the first and last interpolation points are located on the boundary, and the remaining \(n_k - 2\) interpolation points are located inside.

In \eqref{eq:iga_parametric_solution}, the basis functions \(\widehat\varphi_{\mathbf i}\) are taken as the following tensor product B-splines,
\[
\widehat \varphi_{\mathbf i}(\boldsymbol\xi)
=
\prod_{k=1}^{d}N^{(k)}_{i_k,p_k}(\xi^{(k)}),
\qquad
\mathbf i=(i_1,\ldots,i_d)\in\mathcal I.
\]
Here,  \(\boldsymbol{\xi} = (\xi^{(1)}, \ldots, \xi^{(d)})\). The multi-dimensional node vector is
\(
\boldsymbol\Xi
=
\Xi^{(1)}\times\cdots\times\Xi^{(d)},
\) where \(
\Xi^{(k)}
=
\{\xi_1^{(k)},\ldots,\xi_{n_k+p_k+1}^{(k)}\}\).
By selecting one-dimensional Greville collocation points
\(\{\bar{\xi}^{(k)}_{i_k}\}_{i_k=1}^{n_k}\) in each parametric direction and forming their tensor product, we construct the parametric collocation points in \eqref{eq:collocation points},
\[
\bar{\boldsymbol{\xi}}_{\mathbf{i}}
=
\bigl(
\bar{\xi}^{(1)}_{i_1},
\ldots,
\bar{\xi}^{(d)}_{i_d}
\bigr),
\qquad
\mathbf{i}\in\mathcal{I}.
\]

\section{QTT-based IGA Collocation Solver on regular Cartesian domains}

We consider the isogeometric collocation method with QTT representation on regular Cartesian domains.
When the geometry mapping is taken as the identity map, the parameter domain coincides with the physical domain.
In this case, no additional variable coefficients are introduced by the geometry mapping. This case is commonly used as a benchmark for isogeometric collocation methods.

\subsection{Collocation Matrix and Low-rank Operator}

This section assumes that the same number of basis functions \(n\), the same degree \(p\) and the same node vector \(\Xi\) are used in each direction. The discrete solution is expressed as
\begin{equation}
u_h(\boldsymbol{\xi})
=
\sum_{i_1=1}^{n}\cdots\sum_{i_d=1}^{n}
c_{i_1\cdots i_d}
\prod_{k=1}^{d}N_{i_k,p}(\xi^k),
\label{eq:un_expansion}
\end{equation}
where \(c_{i_1\cdots i_d}\) denote the unknown coefficients, and the corresponding tensor is denoted by
\(\boldsymbol{C}=(c_{i_1\cdots i_d})\). To facilitate the QTT representation, let \(n=2^m\), where \(m \in \mathbb N\). Define the one-dimensional collocation matrix
\begin{equation}
\boldsymbol D^{(0)} (i,j)=N_{j,p}(\bar{\xi}_i),
\qquad
\boldsymbol D^{(1)}(i,j)=\frac{d}{d\xi}N_{j,p}(\bar{\xi}_i),
\qquad
\boldsymbol D^{(2)}(i,j)=\frac{d^2}{d\xi^2}N_{j,p}(\bar{\xi}_i),
\label{eq:onedim_mats}
\end{equation}
Here, \(\bar{\xi}_i\) denotes the \(i\)-th Greville abscissa with \(i,j=1,\ldots,n\).
Since the knot vectors are identical in all directions, the matrices \(\boldsymbol D^{(0)},\boldsymbol D^{(1)},\boldsymbol D^{(2)}\) are applicable to every coordinate direction.

Vectorize the coefficient tensor \(\boldsymbol C\) into \(\operatorname{vec}(\boldsymbol C)\in\mathbb R^{N}\) by ordering the last index with the fastest varying rate, where \(N=n^d\), i.e.
\begin{equation}
\operatorname{vec}(\boldsymbol C)\big(\iota(\mathbf i)\big)
=
c_{i_1\cdots i_d},
\qquad
\iota(\mathbf i)
=1+\sum_{k=1}^{d}(i_k-1)n^{d-k}.
\label{eq:vecC}
\end{equation}
For the Greville collocation points in the parameter domain,
\[
\bar{\boldsymbol\xi}_{\mathbf i}
=
(\bar{\xi}_{i_1},\ldots,\bar{\xi}_{i_d}),
\qquad
\mathbf i=(i_1,\ldots,i_d),
\]
it follows from \eqref{eq:un_expansion} and \eqref{eq:onedim_mats} that at the collocation point
\(\bar{\boldsymbol\xi}_{\mathbf i}\),
\begin{equation}
\begin{aligned}
u_h(\bar{\boldsymbol\xi}_{\mathbf i})
&=
\sum_{j_1,\ldots,j_d=1}^{n}
c_{j_1\cdots j_d}
\prod_{k=1}^{d}N_{j_k,p}(\bar{\xi}_{i_k})=
\sum_{j_1,\ldots,j_d=1}^{n}
\Bigl(\prod_{k=1}^{d}\boldsymbol D^{(0)}(i_k,j_k)\Bigr)
c_{j_1\cdots j_d}\\
&=
\Bigl[\bigl(\boldsymbol D^{(0)}\otimes\cdots\otimes \boldsymbol D^{(0)}\bigr)
\operatorname{vec}(\boldsymbol C)\Bigr](\iota(\mathbf i)).
\end{aligned}
\label{eq:un_collocation}
\end{equation}
Let \(\partial_k=\partial/\partial\xi^k\). Since
\(\partial_k\) acts only on the \(k\)-th basis-function factor while the remaining \(d-1\) factors remain unchanged, we have
\begin{equation}
\begin{aligned}
(\partial_k u_h)(\bar{\boldsymbol\xi}_{\mathbf i})
&=
\Bigl[
(\boldsymbol D^{(0)})^{\otimes(k-1)}
\otimes \boldsymbol D^{(1)}
\otimes (\boldsymbol D^{(0)})^{\otimes(d-k)})
\operatorname{vec}(\boldsymbol C)\Bigr](\iota(\mathbf i)),\\
(\partial_k^{2}u_h)(\bar{\boldsymbol\xi}_{\mathbf i})
&=
\Bigl[
(\boldsymbol D^{(0)})^{\otimes(k-1)}
\otimes \boldsymbol D^{(2)}
\otimes (\boldsymbol D^{(0)})^{\otimes(d-k)})
\operatorname{vec}(\boldsymbol C)\Bigr](\iota(\mathbf i)).
\end{aligned}
\label{eq:grad}
\end{equation}

Then the discrete matrix of the Laplace operator
\(\mathcal L=\sum_{k=1}^{d}\partial_k^{2}\) is
\begin{equation}
\boldsymbol L
=
\sum_{k=1}^{d}
(\boldsymbol D^{(0)})^{\otimes(k-1)}
\otimes \boldsymbol D^{(2)}
\otimes( \boldsymbol D^{(0)})^{\otimes(d-k)}.
\label{eq:Lfull}
\end{equation}
Using the method in \cite{kazeev2012}, \eqref{eq:Lfull} can be
recast as
\begin{equation}
\boldsymbol L
=
\begin{bmatrix}
\boldsymbol D^{(0)} & \boldsymbol D^{(2)}
\end{bmatrix}
\bowtie
\begin{bmatrix}
\boldsymbol D^{(0)} & \boldsymbol D^{(2)}\\
\boldsymbol 0 & \boldsymbol D^{(0)}
\end{bmatrix}
\bowtie
\cdots
\bowtie
\begin{bmatrix}
\boldsymbol D^{(0)} & \boldsymbol D^{(2)}\\
\boldsymbol 0 & \boldsymbol D^{(0)}
\end{bmatrix}
\bowtie
\begin{bmatrix}
\boldsymbol D^{(2)}\\
\boldsymbol D^{(0)}
\end{bmatrix}.
\label{eq:Lblock}
\end{equation}
Assume that the one-dimensional matrices have QTT representations,
\[
\boldsymbol D^{(0)}
=
\boldsymbol U_1\bowtie \boldsymbol U_2\bowtie\cdots\bowtie \boldsymbol U_m,
\qquad
\boldsymbol D^{(2)}
=
\boldsymbol V_1\bowtie \boldsymbol V_2\bowtie\cdots\bowtie \boldsymbol V_m,
\]
Here, the construction of \(\boldsymbol U_i\) and \(\boldsymbol V_i\) is deferred to the next subsection.  Arrange these  core matrice as follows,
\begin{align}
\begin{bmatrix}
\boldsymbol D^{(0)} & \boldsymbol D^{(2)}\\
\boldsymbol 0 & \boldsymbol D^{(0)}
\end{bmatrix}
&=
\left[
\begin{array}{ccc}
\boldsymbol U_1 & \boldsymbol V_1 & \\
& & \boldsymbol U_1
\end{array}
\right]
\bowtie
\left[
\begin{array}{ccc}
\boldsymbol U_2 & & \\
& \boldsymbol V_2 & \\
& & \boldsymbol U_2
\end{array}
\right]
\bowtie
\cdots
\bowtie
\left[
\begin{array}{ccc}
\boldsymbol U_m & \\
& \boldsymbol V_m \\
& \boldsymbol U_m
\end{array}
\right],
\label{eq:block_qtt_mid}
\\[1em]
\left[
\begin{array}{cc}
\boldsymbol D^{(0)} & \boldsymbol D^{(2)}
\end{array}
\right]
&=
\left[
\begin{array}{cc}
\boldsymbol U_1 & \boldsymbol V_1
\end{array}
\right]
\bowtie
\left[
\begin{array}{cc}
\boldsymbol U_2 & \\
& \boldsymbol V_2
\end{array}
\right]
\bowtie
\cdots
\bowtie
\left[
\begin{array}{cc}
\boldsymbol U_m & \\
& \boldsymbol V_m
\end{array}
\right],
\label{eq:block_qtt_row}
\\[1em]
\left[
\begin{array}{c}
\boldsymbol D^{(2)}\\
\boldsymbol D^{(0)}
\end{array}
\right]
&=
\left[
\begin{array}{cc}
\boldsymbol V_1 & \\
& \boldsymbol U_1
\end{array}
\right]
\bowtie
\cdots
\bowtie
\left[
\begin{array}{cc}
\boldsymbol V_{m-1} & \\
& \boldsymbol U_{m-1}
\end{array}
\right]
\bowtie
\left[
\begin{array}{cc}
\boldsymbol V_m \\
\boldsymbol U_m
\end{array}
\right],
\label{eq:block_qtt_col}
\end{align}
Substituting these equations into \eqref{eq:Lblock} yields the QTT-matrix
representation of \(\boldsymbol L\). 
If different \(n_k, p_k\) and node vectors are used in each direction, the corresponding QTT-matrix representation can still be obtained by the similar  method described in~\cite{kazeev2012}.

\subsection{QTT Representation of One-Dimensional Collocation Matrices}

This section discusses the QTT representations of the matrices
\(\boldsymbol D^{(0)}\), \(\boldsymbol{D}^{(1)}\), and \(\boldsymbol{D}^{(2)}\).
These three matrices share a common banded structure and interior translation invariance, which can be exploited to construct their QTT decompositions and to derive rank estimates.

\subsubsection{Interior Translation Invariance}

According to the definition of the Greville abscissa and the support set of B-splines of degree \(p\), the matrices
\(\boldsymbol D^{(0)}\), \(\boldsymbol{D}^{(1)}\), and \(\boldsymbol{D}^{(2)}\)
are banded, with semi-bandwidth at most \(p\).
For \(k=-p,\ldots,p\), define the entries on the \(k\)-th band by
\[
\boldsymbol{d}_{k}^{(s)}(i)
=
\boldsymbol{A}^{(s)}(i,i+k),
\qquad
i\in I_k,
\]
where \(I_k\) denotes the set of row indices contained in the \(k\)-th band, namely
\[
I_k
=
\{i\in\{1,\ldots,n\}:1\le i+k\le n\}
=
[\max(1,1-k),\min(n,n-k)].
\]
\begin{theorem}
Under the above uniform node setting, let
\[
P_k=\bigl[\max(p,\,p+1-k),\ \min(n-p+1,\,n-p-k)\bigr].
\]
If \(P_k \ne \varnothing\), then there exists a constant \(c_k^{(s)}\) such that
\(\boldsymbol d_{k}^{(s)}(i) = c_k^{(s)}\) holds for all \(i \in P_k\) and \(s = 0, 1, 2\). Furthermore,
\(|I_k\setminus P_k|\le 2p\).
\end{theorem}

\begin{proof}
Let \(h = \dfrac{1}{n-p}\), and the uniform node vector is
\begin{equation*}
\xi_1=\cdots=\xi_{p+1}=0,\,
\xi_{p+1+j}=jh\ \ (j=1,\ldots,n-p-1),\,
\xi_{n+1}=\cdots=\xi_{n+p+1}=1.
\end{equation*}
 \(\xi_r = (r - p - 1)h\) holds for \(p + 1 \le r \le n + 1\).

We first prove that for \(p+1\le j\le n-p'\) and $p'=0,1,\ldots,p$,
\begin{equation}
    N_{j,p'}(x)=N_{p+1,p'}\bigl(x-(j-p-1)h\bigr)\label{eq:Invariance}.
\end{equation}
The proof proceeds by induction. For \(p'=0\), both sides correspond to indicator functions over \([\xi_j,\xi_{j+1})\) and \([\xi_{p+1},\xi_{p+2})\), respectively.  Since \(p+1\le j\le n\), the interval \([\xi_j,\xi_{j+1})\) is contained within \([\xi_{p+1},\xi_{n+1}]\), so we immediately obtain \eqref{eq:Invariance} for \(p=0\).
Assume the statement holds for \(p'=p-1\), for \(p+1\le j\le n-p+1\),
\[
N_{j,p-1}(x)=N_{p+1,p-1}\bigl(x-(j-p-1)h\bigr).
\]
Let \(x'=x-(j-p-1)h\). For \(p'=p\) and $p+1\le j\le n-p$,
\begin{equation}
\begin{aligned}
N_{j,p}(x)
&=
\frac{x-\xi_j}{\xi_{j+p}-\xi_j}
N_{j,p-1}(x)
+
\frac{\xi_{j+p+1}-x}{\xi_{j+p+1}-\xi_{j+1}}
N_{j+1,p-1}(x)\\
&=
\frac{x'-\xi_{p+1}}{\xi_{2p+1}-\xi_{p+1}}
N_{p+1,p-1}(x')
+
\frac{\xi_{2p+2}-x'}{\xi_{2p+2}-\xi_{p+2}}
N_{p+2,p-1}(x')
=
N_{p+1,p}(x').
\end{aligned}
\end{equation}
Since \(p+1\le j\le n-p\), we have \(\xi_{p+1}=0\), \(\xi_j=(j-p-1)h\), \(\xi_{j+p+1}=jh\), and \(\xi_{2p+2}=(p+1)h\).Consequently, we obtain \(x-\xi_j=x'-\xi_{p+1}\) and \(\xi_{j+p+1}-x=\xi_{2p+2}-x'\). Moreover,
\[
\xi_{j+p}-\xi_j=\xi_{j+p+1}-\xi_{j+1}=\xi_{2p+1}-\xi_{p+1}=\xi_{2p+2}-\xi_{p+2}=ph,
\]
which implies that the two pairs of coefficients in the second equality are identical. 
Hence,
\begin{equation}
\frac{d^s}{dx^s}N_{j,p}(x)=\frac{d^s}{dx^s}N_{p+1,p}\bigl(x-(j-p-1)h\bigr),
\quad s = 1, 2.
\label{eq:shift}
\end{equation}

The \(k\)-th band is associated with the basis function \(N_{i+k,p}\). Setting \(j=i+k\), then
\begin{align}
    N_{i+k,p}(x)
&=
N_{p+1,p}\bigl(x-(i+k-p-1)h\bigr),\\
\,
\frac{d^s}{dx^s}N_{i+k,p}(x)
&=
\frac{d^s}{dx^s}N_{p+1,p}\bigl(x-(i+k-p-1)h\bigr)
\label{eq:shiftk}
\end{align}
hold if and only if
\begin{equation}
p+1\le i+k\le n-p,
\label{eq:rangek}
\end{equation}
for \(s=0,1,2\).

On the other hand, according to the definition of Greville abscissa,
\begin{equation}
\bar{\xi}_i-\bar{\xi}_{p}=(i-p)h,\qquad p\le i\le n-p+1 .
\label{eq:grevrange}
\end{equation}

Equation \eqref{eq:rangek} and Equation \eqref{eq:grevrange} being simultaneously valid is equivalent to \(i \in P_k\). For such \(i\), setting \(x=\bar{\xi}_i\) in \eqref{eq:shiftk} and using \eqref{eq:grevrange} to replace \(\bar{\xi}_i\) by \(\bar{\xi}_{p}+(i-p)h\), we obtain
\[
\begin{aligned}
\boldsymbol{d}_{k}^{(s)}(i)
&=\frac{d^s}{dx^s}N_{i+k,p}(\bar{\xi}_i)
=\frac{d^s}{dx^s}N_{p+1,p}\bigl(\bar{\xi}_i-(i+k-p-1)h\bigr)\\
&=\frac{d^s}{dx^s}N_{p+1,p}\Bigl(\bar{\xi}_{p}+(i-p)h-(i+k-p-1)h\Bigr)
=\frac{d^s}{dx^s}N_{p+1,p}\bigl(\bar{\xi}_{p}-(k-1)h\bigr).
\end{aligned}
\]
The right-hand side is independent of \(i\). Hence, \(\boldsymbol{d}_{k}^{(s)}\) is constant on \(P_k\), and we denote this constant by \(c_k^{(s)}\).

Next, we need to discuss \(I_k \setminus P_k\). The three cases \(1\le k\le p\), \(k=0\), and \(-p\le k\le-1\) yield, respectively,
\begin{align*}
    I_k&=[1,n-k],\, P_k=[p,n-p-k];\,\,\\
I_0&=[1,n],\, P_0=[p+1,n-p];\,\,\\
I_k&=[1-k,n],\, P_k=[p+1-k,n-p+1].
\end{align*}
In all three cases, \(P_k\) is shorter than \(I_k\) by \(p\) rows at the left and right ends. Consequently, whenever \(P_k\ne\varnothing\), we have \(|I_k\setminus P_k|\leq2p\); if \(P_k=\varnothing\), then the inequalities
\(n<2p+k\), \(n<2p+1\), and \(n<2p-k\) give \(|I_k|\le 2p\). Therefore, \(|I_k\setminus P_k|\le 2p\). This completes the proof.
\end{proof}

The interval \(P_k\) is the smallest interval that ensures the elements in the $k$-th band remain constant. The actual constant interval may be larger. 

\subsubsection{The algorithm for constructing QTT}

Define \( \boldsymbol{\bar{d}}_k^{(s)} \in \mathbb R^n \), where \( k = -p, \ldots, p, \, s = 1, 2, 3 \), such that
\[
\boldsymbol{\bar{d}}_{k}^{(s)}(i)
=
\begin{cases}
\boldsymbol{d}_{k}^{(s)}(i), & i\in I_k,\\
c_k^{(s)}, & i\in[1,n]\setminus I_k.
\end{cases}
\]
The non-cyclic shift matrix $\boldsymbol S_k \in \mathbb R^{n\times n}$ is defined as
\[
\boldsymbol S_k(i,j)
=
\begin{cases} 
1, & j = i + k \ \text{and} \ 1 \le j \le n, \\
0, & \text{otherwise}.
\end{cases}
\]
Therefore,
\begin{equation}
\boldsymbol A^{(s)}=\sum_{k=-p}^{p}\operatorname{diag}(\boldsymbol{\bar{d}}_k^{(s)})\boldsymbol S_k.
\label{eq:decom}
\end{equation}

When constructing the QTT matrix, only a small number of matrix elements need to be accessed. First, select one row from the \(k\)-th band of \(\boldsymbol A^{(s)}\) within the interval \(P_k\) to obtain \(c_k^{(s)}\). Then, access each element in \(I_k \setminus P_k\) one by one and compare it with \(c_k^{(s)}\). For \(j \in I_k \setminus P_k\), define \(\boldsymbol \rho_{k}^{(s)}(j) = \boldsymbol d_{k}^{(s)}(j) - c_k^{(s)}\). Keep the elements that satisfy the following conditions,
\[
|\boldsymbol \rho_{k}^{(s)}(j)| > \varepsilon \max\{|c_k^{(s)}|, |\boldsymbol d_{k}^{(s)}(j)|, 1\}.
\]
Here, \(\varepsilon\) represents the small tolerance. The set of retained indicators is denoted as
\[
\mathcal{E} = \{j \in I_k \setminus P_k : |\boldsymbol \rho_{k}^{(s)}(j)| > \varepsilon \max\{|c_k^{(s)}|, |\boldsymbol d_{k}^{(s)}(j)|, 1|\}
\] 
Therefore,
\[
\boldsymbol{\bar{d}}_k^{(s)}
=
c_k^{(s)}\mathbf 1 + \sum_{j\in\mathcal E}\boldsymbol\rho_{k}^{(s)}(j)\boldsymbol e_j.
\]
Here, $\boldsymbol e_j \in \mathbb R^n$ denotes the $j$-th standard basis vector.

We next construct the QTT matrix corresponding to \(\operatorname{diag}(\boldsymbol{\bar{d}}_k^{(s)})\). First,
\(\mathbf 1 = \bigotimes_{\ell=1}^{m}\begin{pmatrix}1&1\end{pmatrix}\). Second, for any \(j \in \{1, \ldots, n\}\),  let \(
j=1+\sum_{k=1}^{m}\beta_k2^{m-k},
 i_k\in\{0,1\}.
\)
 Denote 
\(\boldsymbol e_0 = (1, 0)^{\mathrm T}\) and \(\boldsymbol e_1 = (0, 1)^{\mathrm T}\). Then,
\(\boldsymbol e_j = \boldsymbol e_{\beta_1} \otimes \cdots \otimes \boldsymbol e_{\beta_m}\). Therefore, both \(\mathbf 1\) and \(\boldsymbol e_j\) are QTT vectors of rank 1. Adding \(c_k^{(s)}\mathbf 1\) entrywise to each \(\boldsymbol\rho_{k}^{(s)}(j)\,\boldsymbol e_j\) and performing TT rounding at each step yields a QTT vector of \(\boldsymbol{\bar{d}}_k^{(s)}\), denoted by \(\boldsymbol{\mathbf X}_k^{(s)}\), whose \(\ell\)-th core is denoted by
\(\boldsymbol X_{\ell,k}^{(s)} \in \mathbb R^{r_{\ell-1}\times2\times r_\ell}\). Its rank satisfies
\[
\rank_{\mathrm{QTT}}(\boldsymbol{\mathbf X}_k^{(s)})\le1+|\mathcal E|\le1+2p.
\]
The cores of the diagonal QTT matrix obtained from this TT vector are
\begin{equation}
\boldsymbol W_{\ell,k}^{(s)}(\alpha,i,j,\gamma)
=
\boldsymbol X_{\ell,k}^{(s)}(\alpha,i,\gamma)\delta_{ij},   
\label{eq:diag_QTT}
\end{equation}
where $\delta_{ij}$ denotes the Kronecker delta. Hence,
$\operatorname{diag}(\boldsymbol{\bar{d}}_k^{(s)}) = \boldsymbol W_{1,k}^{(s)} \bowtie \cdots \bowtie \boldsymbol W_{m,k}^{(s)}$.

Then, according to \eqref{eq:decom}, we construct the QTT matrix of \(\boldsymbol A^{(s)}\). The QTT rank of the non-cyclic shift matrix \(\boldsymbol S_k\) is at most \(2\) \cite{kazeev2013}. By multiplying the QTT representation of \(\operatorname{diag}(\boldsymbol{\bar{d}}_k^{(s)})\) with that of \(\boldsymbol S_k\), we obtain
\[
\operatorname{rank}_{\mathrm{QTT}}\bigl(\operatorname{diag}(\boldsymbol{\bar{d}}_k^{(s)})\boldsymbol S_k\bigr)
\le 2(1+2p)=O(p).
\]
Summing over \(k=-p,\ldots,p\) and performing TT rounding at each step yields the QTT matrix of \(\boldsymbol A^{(s)}\). Consequently,
\[
\operatorname{rank}_{\mathrm{QTT}}\bigl(\boldsymbol A^{(s)}\bigr)=O(p^2).
\]
In the above construction, for each band it suffices to access \(p\) rows at each of the two ends and only one row in the interior. Since there are \(2p+1\) bands in total, the total number of accessed entries is \(O(p^2)\). The storage of the corresponding one-dimensional QTT matrix is \(O(p^4\log n)\).

\subsection{QTT Representation of Boundary conditions}

We use the two-dimensional Poisson problem as a representative example to
explain how Dirichlet and Neumann boundary conditions enter the tensor
formulation. Let \(\hat{\Omega}=(0,1)^{2}\) and consider
\begin{equation}
    -\Delta u(\xi^{1},\xi^{2})=f(\xi^{1},\xi^{2})
    \qquad\text{in }\hat{\Omega}.
    \label{eq:PDE}
\end{equation}
Dirichlet conditions are imposed on the two boundaries in the first
parameter direction,
\begin{equation}
    u(0,\xi^{2})=h_1(\xi^{2}),
    \qquad
    u(1,\xi^{2})=h_2(\xi^{2}),
    \label{eq:BCdir}
\end{equation}
 while Neumann conditions are imposed on the remaining
two boundaries,
\begin{equation}
    -\partial_{2}u(\xi^{1},0)=d_1(\xi^{1}),
    \qquad
    \partial_{2}u(\xi^{1},1)=d_2(\xi^{1}),
    \label{eq:BCneu}
\end{equation}
where \(\partial_2=\partial/\partial\xi^2\). The functions
\(h_1,h_2,d_1,d_2\) are prescribed boundary data and are assumed to
satisfy the corner compatibility conditions required by the continuous
problem. 
The discrete solution is represented as
\(
    u_h(\xi^{1},\xi^{2})
    =
    \sum_{i=1}^{n}\sum_{j=1}^{n}
    c_{ij}N_{i,p}(\xi^{1})N_{j,p}(\xi^{2}),
\)
where the coefficient matrix is
\(
    \boldsymbol C=(c_{ij}).
\)

\subsubsection{Elimination of boundary unknowns}

Because the open uniform knot vectors interpolate the endpoints, the
univariate B-splines satisfy \(   N_{i,p}(0)=\delta_{i,1},
    N_{i,p}(1)=\delta_{i,n}.
    \label{eq:endpointinterp}
\)
Consequently,
\begin{equation}
     u_h(0,\xi^{2})
    =
    \sum_{j=1}^{n}c_{1j}N_{j,p}(\xi^{2}),\qquad
    u_h(1,\xi^{2})
    =
    \sum_{j=1}^{n}c_{nj}N_{j,p}(\xi^{2}),
    \label{eq:dy1}
\end{equation}
Substituting \eqref{eq:dy1} into the Dirichlet condition \eqref{eq:BCdir} yields
\begin{equation*}
    \sum_{j=1}^{n}\boldsymbol D^{(0)}(i,j)c_{1,j}
    =
    h_1(\bar{\xi}_i),
    \qquad
    \sum_{j=1}^{n}\boldsymbol D^{(0)}(i,j)c_{n,j}
    =
    h_2(\bar{\xi}_i),
    \qquad
    i=1,\ldots,n.
\end{equation*}
Define 
\begin{equation}
    \boldsymbol{\alpha}^{(1)}
    =
    (\boldsymbol D^{(0)})^{-1}\boldsymbol h_1,
    \qquad
    \boldsymbol{\alpha}^{(2)}
    =
    (\boldsymbol D^{(0)})^{-1}\boldsymbol h_2.
    \label{eq:alpha}
\end{equation}
where 
\(
    \boldsymbol h_1=    \bigl(h_1(\bar{\xi}_1),\ldots,h_1(\bar{\xi}_n)\bigr)^{\mathrm T},
    \boldsymbol h_2=    \bigl(h_2(\bar{\xi}_1),\ldots,h_2(\bar{\xi}_n)\bigr)^{\mathrm T}.
\)
Hence,
\begin{equation}
    c_{1,i}=\alpha_i^{(1)},
    \qquad
    c_{n,i}=\alpha_i^{(2)},
    \qquad
    i=1,\ldots,n.
    \label{eq:Cdir}
\end{equation}
where \(\alpha_i^{(1)}\) and \(\alpha_i^{(2)}\) are the \(i\)-th entries
of \(\boldsymbol{\alpha}^{(1)}\) and \(\boldsymbol{\alpha}^{(2)}\),
respectively.

For the Neumann boundaries,  the endpoint
derivatives of the B-splines satisfy
\[
    N_{1,p}'(0)=-\frac{p}{h},
    \qquad
    N_{2,p}'(0)=\frac{p}{h},
    \qquad
    N_{i,p}'(0)=0,\quad i\geq3,
\]
and
\[
    N_{n-1,p}'(1)=-\frac{p}{h},
    \qquad
    N_{n,p}'(1)=\frac{p}{h},
    \qquad
    N_{i,p}'(1)=0,\quad i\leq n-2,
\]
where \(h = 1/(n - p)\). Therefore,
\begin{equation}
    \partial_{2}u_h(\xi^{1},0)
    =
    \frac{p}{h}
    \sum_{k=1}^{n}
    \bigl(c_{k,2}-c_{k,1}\bigr)N_{k,p}(\xi^{1}),
   \,
    \partial_{2}u_h(\xi^{1},1)
    =
    \frac{p}{h}
    \sum_{k=1}^{n}
    \bigl(c_{k,n}-c_{k,n-1}\bigr)N_{k,p}(\xi^{1}).
    \label{eq:dy}
\end{equation}
Substituting \eqref{eq:dy} into the Neumann condition \eqref{eq:BCneu} yields
\begin{equation*}
    \frac{p}{h}
    \sum_{k=1}^{n}\boldsymbol D^{(0)}(i,k)
    \bigl(c_{k,1}-c_{k,2}\bigr)
    =
    d_1(\bar{\xi}_i),
    \,
    \frac{p}{h}
    \sum_{k=1}^{n}\boldsymbol D^{(0)}(i,k)
    \bigl(c_{k,n}-c_{k,n-1}\bigr)
    =
    d_2(\bar{\xi}_i),
\end{equation*}
where $i=1,\ldots,n.$
Define
\begin{equation}
    \boldsymbol{\beta}^{(1)}
    =
    \frac{h}{p}\,(\boldsymbol D^{(0)})^{-1}\boldsymbol d_1,
    \qquad
    \boldsymbol{\beta}^{(2)}
    =
    \frac{h}{p}\,(\boldsymbol D^{(0)})^{-1}\boldsymbol d_2.
    \label{eq:beta}
\end{equation}
where \(
    \boldsymbol d_1=
\bigl(d_1(\bar{\xi}_1),\ldots,d_1(\bar{\xi}_n)\bigr)^{\mathrm T},
    \boldsymbol d_2=
\bigl(d_2(\bar{\xi}_1),\ldots,d_2(\bar{\xi}_n)\bigr)^{\mathrm T},
\)
Thus,  
\begin{equation}
    c_{i,1}-c_{i,2}=\beta_i^{(1)},
    \qquad
    c_{i,n}-c_{i,n-1}=\beta_i^{(2)},
    \qquad
    i=1,\ldots,n,
    \label{eq:neumannrow1}
\end{equation}
where \(\beta_i^{(1)}\) and \(\beta_i^{(2)}\) are the \(i\)-th entries
of \(\boldsymbol{\beta}^{(1)}\) and \(\boldsymbol{\beta}^{(2)}\),
respectively.
Taking into account the compatibility conditions at the corners, the Neumann boundary conditions only need to be imposed at the 
\(n-2\)  interior collocation points. Therefore, we introduce a zero-padded boundary correction vector, 
\(
    \widetilde{\boldsymbol{\beta}}^{(1)}
    =
    (\widetilde{\beta}_1^{(1)},\ldots,
    \widetilde{\beta}_n^{(1)})^{\mathrm T},
    \widetilde{\boldsymbol{\beta}}^{(2)}
    =
    (\widetilde{\beta}_1^{(2)},\ldots,
    \widetilde{\beta}_n^{(2)})^{\mathrm T},
    \label{eq:betatildevec}
\)
with entries
\begin{equation}
    \widetilde{\beta}_i^{(1)}
    =
    \begin{cases}
        \beta_i^{(1)}, & i=2,\ldots,n-1,\\[1mm]
        0, & i=1\ \text{or}\ i=n,
    \end{cases}
    \qquad
    \widetilde{\beta}_i^{(2)}
    =
    \begin{cases}
        \beta_i^{(2)}, & i=2,\ldots,n-1,\\[1mm]
        0, & i=1\ \text{or}\ i=n.
    \end{cases}
    \label{eq:betatilde}
\end{equation}

We next introduce an interior coefficient matrix
\(\boldsymbol C^{(\mathrm{int})}=(c^{(\mathrm{int})}_{i,j})\) , with all coefficients associated with the boundary layers set to zero, so that the dimension remains an integer power of 2,
\begin{equation}
    c^{(\mathrm{int})}_{i,j}
    =
    \begin{cases}
        c_{i,j}, & 2\leq i,j\leq n-1,\\[1mm]
        0, & \text{otherwise}.
    \end{cases}
    \label{eq:Ydef}
\end{equation}
Define the one-dimensional matrices,
\begin{equation}
    \boldsymbol E_D
    =
    \operatorname{diag}(0,1,\ldots,1,0) \in \mathbb R^{n\times n},
    \qquad
    \boldsymbol E_N
    =
    \boldsymbol E_D+\boldsymbol e_1\boldsymbol e_2^{\mathrm T}
    +\boldsymbol e_n\boldsymbol e_{n-1}^{\mathrm T},
    \label{eq:DI}
\end{equation} and \(
    \boldsymbol E=\boldsymbol E_D\otimes\boldsymbol E_N.
\)
Using \eqref{eq:Cdir} and \eqref{eq:betatilde}, the four boundary layers
are collected in the lifting matrix
\(\boldsymbol C^{(\mathrm{bd})}=(c^{(\mathrm{bd})}_{i,j})\), with entries
\begin{equation}
    \begin{aligned}
    c^{(\mathrm{bd})}_{i,j}
    ={}&
    \alpha_j^{(1)}\delta_{i,1}
    +
    \alpha_j^{(2)}\delta_{i,n}
    +
    \widetilde{\beta}_i^{(1)}\delta_{j,1}
    +
    \widetilde{\beta}_i^{(2)}\delta_{j,n}.
    \end{aligned}
    \label{eq:C0}
\end{equation}
Combining \eqref{eq:Cdir} and \eqref{eq:neumannrow1}, it can be verified that the complete coefficient vector is expressed as
\begin{equation}
    \operatorname{vec}(\boldsymbol C)
    =
    \boldsymbol E\operatorname{vec}(\boldsymbol C^{(\mathrm{int})})
    +
    \operatorname{vec}(\boldsymbol C^{(\mathrm{bd})}).
    \label{eq:Cdecomp}
\end{equation}
If we can obtain $\boldsymbol C^{(\mathrm{int})}$ which includes the internal coefficients, then by using \eqref{eq:Cdecomp}, we can recover the complete coefficients $\boldsymbol C$.

\subsubsection{QTT representation and solution of the reduced system}

The discrete equation that does not involve boundary condition handling is $\boldsymbol L\operatorname{vec}(\boldsymbol C) = \operatorname{vec}(\boldsymbol F)$.
where \(\boldsymbol L=\boldsymbol D^{(2)}\otimes\boldsymbol D^{(0)} +\boldsymbol D^{(0)} \otimes \boldsymbol D^{(2)} \),
\(\boldsymbol F=(f_{i,j})\) is the right-hand side, and
\(f_{i,j}=f(\bar{\xi}_i,\bar{\xi}_j)\).
Next, the equations associated with the boundary collocation points are set to trivial identities $0=0$, so that they are effectively removed while the total number of equations remains an integer power of 2. Define
\(
    \boldsymbol M
    =
    \boldsymbol E_D\otimes\boldsymbol E_D
    \label{eq:twodimask}
\) and then
the interior equations can be expressed as
\begin{equation}
    \boldsymbol M\boldsymbol L\operatorname{vec}(\boldsymbol C)
    =
    \boldsymbol M\operatorname{vec}(\boldsymbol F),
    \label{eq:maskedsystem}
\end{equation}
Substituting
\eqref{eq:Cdecomp} into \eqref{eq:maskedsystem} gives
\begin{equation*}
    \boldsymbol M\boldsymbol L\boldsymbol E
    \operatorname{vec}(\boldsymbol C^{(\mathrm{int})})
    =
    \boldsymbol M
    \bigl(
    \operatorname{vec}(\boldsymbol F)
    -
    \boldsymbol L\operatorname{vec}(\boldsymbol C^{(\mathrm{bd})})
    \bigr).
\end{equation*}
Set
\begin{equation*}
    \widehat{\boldsymbol L}
    =
    \boldsymbol M\boldsymbol L\boldsymbol E,
    \qquad
    \widehat{\boldsymbol f}
    =
    \boldsymbol M
    \bigl(
    \operatorname{vec}(\boldsymbol F)
    -
    \boldsymbol L\operatorname{vec}(\boldsymbol C^{(\mathrm{bd})})
    \bigr),
\end{equation*}
so that
\begin{equation}
    \widehat{\boldsymbol L}\operatorname{vec}(\boldsymbol C^{(\mathrm{int})})
    =
    \widehat{\boldsymbol f}.
    \label{eq:reduced_masked}
\end{equation}
 All rows
associated with padded boundary slots are zero. Consequently,
\(\widehat{\boldsymbol L}\) is singular. We therefore add the
term \(\boldsymbol I-\boldsymbol M\) and solve the equivalent
nonsingular system,
\begin{equation}
    \widehat{\boldsymbol L}_{\mathrm{solve}}
    \operatorname{vec}({\boldsymbol C^{(\mathrm{int})}})
    =
    \widehat{\boldsymbol f},
    \qquad
    \widehat{\boldsymbol L}_{\mathrm{solve}}
    =
    \widehat{\boldsymbol L}+\boldsymbol I-\boldsymbol M.
    \label{eq:Ahatsys}
\end{equation}
In the boundary blocks, the corresponding rows of \(\widehat{\boldsymbol L}_{\mathrm{solve}}\) are identity rows and the entries of the right-hand side \(\widehat{\boldsymbol f}\) are zero at these positions, so the values of the solution \({\boldsymbol C^{(\mathrm{int})}}\) at the corresponding positions are zero. In the interior blocks, the values of the solution \({\boldsymbol C^{(\mathrm{int})}}\) from \eqref{eq:Ahatsys} coincide with those of the solution \(\boldsymbol C^{(\mathrm{int})}\) from \eqref{eq:reduced_masked}.

Next, we will explain the QTT representation involved in \eqref{eq:Ahatsys}. Since the QTT representations of \(\boldsymbol E_D=\boldsymbol I-\boldsymbol e_1 \boldsymbol e_1^{\mathrm T}-\boldsymbol e_n \boldsymbol e_n^{\mathrm T}\) and \(\boldsymbol E_N =
    \boldsymbol E_D+\boldsymbol e_1\boldsymbol e_2^{\mathrm T}
    +\boldsymbol e_n\boldsymbol e_{n-1}^{\mathrm T}\) are both easy to construct, the QTT representations of \(\boldsymbol M\) and \(\boldsymbol E\) follow. Meanwhile, the QTT matrix representation of \(\boldsymbol L\) has been described in detail in the previous two sections. Combining these ingredients, the QTT representation of \(\widehat{\boldsymbol L}_{\mathrm{solve}}\) is readily obtained.
 On the right-hand side, the vector \(\operatorname{vec}(\boldsymbol F)\) can be constructed into the QTT representation through the TT-cross approximation. By \eqref{eq:C0}, the lifting item \(\boldsymbol C^{(\mathrm{bd})}\) is the sum of four rank-one terms,
\begin{equation}
    \boldsymbol C^{(\mathrm{bd})}
    =
    \boldsymbol e_1\otimes\boldsymbol{\alpha}^{(1)}
    +
    \boldsymbol e_n\otimes\boldsymbol{\alpha}^{(2)}
    +
    \widetilde{\boldsymbol{\beta}}^{(1)}\otimes\boldsymbol e_1
    +
    \widetilde{\boldsymbol{\beta}}^{(2)}\otimes\boldsymbol e_n.
    \label{eq:C0sum}
\end{equation}
The QTT representations of the  boundary coefficient vectors \(\boldsymbol{\alpha}^{(1)},\boldsymbol{\alpha}^{(1)},\widetilde{\boldsymbol{\beta}}^{(1)},\widetilde{\boldsymbol{\beta}}^{(1)}\) can be obtained by solving \eqref{eq:alpha} and \eqref{eq:beta} with the TT-solver, where the QTT matrix of \(\boldsymbol D^{(0)}\) has already been constructed in the previous section. Finally, solving \eqref{eq:Ahatsys} with the TT-solver yields the QTT representation of \({\boldsymbol C^{(\mathrm{int})}}\). The QTT representation of the complete coefficient vector is reconstructed from \eqref{eq:Cdecomp}.

The above handling of the boundary conditions is presented in two dimensions only for notational convenience. The $d$-dimensional case can be handled in a completely analogous manner.


\section{QTT-based IGA Collocation Solver for Poisson Problems on 3D nontrivial geometries}
\label{sec:mapped}

In the preceding section, the parameter domain was assumed to coincide with the physical domain, with the geometry mapping taken as the identity map. This section removes this restriction. Let the physical domain \(\Omega\subset\mathbb R^3\) be obtained by a smooth and invertible mapping \(\boldsymbol G:[0,1]^3\to\overline\Omega\).

We consider the 3D Poisson problem as in \cite{tran2026},
\begin{equation}
-\Delta u(\boldsymbol x)=f(\boldsymbol x)\quad\text{in }\Omega,
\label{eq:poisson}
\end{equation}
where \(u(\boldsymbol x)\) denotes the scalar field to be solved for and \(f(\boldsymbol x)\) is a known function.
The boundary \(\partial\Omega\) is 
\begin{equation}
\begin{cases}
u(\boldsymbol x)=g_D(\boldsymbol x), & \boldsymbol x\in\Gamma_D\subset\partial\Omega\quad(\text{Dirichlet condition}),\\[4pt]
\mathbf{n}\cdot\nabla u(\boldsymbol x)=g_N(\boldsymbol x), & \boldsymbol x\in\Gamma_N=\partial\Omega\setminus\Gamma_D\quad(\text{Neumann condition}),
\end{cases}
\label{eq:mappde}
\end{equation}
where \(\mathbf{n}\) stands for the outward-pointing unit normal vector to the boundary, \(g_D(\boldsymbol x)\) is the prescribed boundary data on the Dirichlet boundary \(\Gamma_D\), and \(g_N(\boldsymbol x)\) denotes the prescribed  Neumann traction on  \(\Gamma_N\).

\subsection{Parameter Representation of Poisson Equation}

Let 
\(
\widehat u(\boldsymbol\xi)=u\bigl(\boldsymbol G(\boldsymbol\xi)\bigr)
\) and define \(\boldsymbol G(\boldsymbol \xi)=(G_1(\boldsymbol \xi),G_2(\boldsymbol \xi),G_3(\boldsymbol \xi))^{\mathrm T}\). 
In this section, we reformulate \eqref{eq:mappde} equivalently as a problem for \(\widehat u\) defined over \([0,1]^3\), where all geometric contributions are contained in several variable coefficients.

Let
\begin{equation}
J_{sr}(\boldsymbol{\xi})
=\frac{\partial G_s}{\partial \xi^r}(\boldsymbol{\xi}),
\qquad
\boldsymbol J(\boldsymbol{\xi})
=\bigl(J_{sr}(\boldsymbol{\xi})\bigr)_{s,r=1}^{3},
\label{eq:jacdef}
\end{equation}
where \(\boldsymbol J(\boldsymbol{\xi})\) is the Jacobian matrix, and let \(J^{-1}_{sr}(\boldsymbol{\xi})\) stand for the \((s,r)\)-th entry of \(\bigl(\boldsymbol J(\boldsymbol{\xi})\bigr)^{-1}\).
We then have
\begin{equation}
\frac{\partial u}{\partial x_r}\bigl(\boldsymbol G(\boldsymbol\xi)\bigr)
=\sum_{s=1}^3 J^{-1}_{sr} (\boldsymbol\xi)\,
\frac{\partial\widehat u}{\partial\xi^s}(\boldsymbol\xi),
\qquad r=1,2,3 .
\label{eq:firstdercomp}
\end{equation}

From \eqref{eq:firstdercomp}, the first-order differential operator in the physical direction can be expressed as
\begin{equation}
\frac{\partial}{\partial x_r}
=\sum_{q=1}^3 J^{-1}_{qr}(\boldsymbol\xi)\,
\frac{\partial}{\partial\xi^q},
\qquad r=1,2,3 .
\label{eq:physderparam}
\end{equation}
Applying this operator to \eqref{eq:firstdercomp} and  yielding
\begin{equation}
\frac{\partial^2 u}{\partial x_r^2}\bigl(\boldsymbol G(\boldsymbol\xi)\bigr)
=\sum_{s,t=1}^3 J^{-1}_{sr}(\boldsymbol\xi)J^{-1}_{tr}(\boldsymbol\xi)\,
\frac{\partial^2\widehat u}{\partial\xi^t\partial\xi^s}
+\sum_{s=1}^3\Bigl(\sum_{q=1}^3  J^{-1}_{qr}(\boldsymbol\xi)\,
\partial_q\big( J^{-1}_{sr}(\boldsymbol\xi)\big)\Bigr)\frac{\partial\widehat u}{\partial\xi^s}.
\label{eq:hessder}
\end{equation}

Summing over \(r = 1, 2, 3\), we obtain the form of the Laplace operator in the parameter coordinates,
\begin{equation}
\Delta u\bigl(\boldsymbol G(\boldsymbol\xi)\bigr)
=\sum_{s,t=1}^3 a_{st}(\boldsymbol\xi)\,
\frac{\partial^2\widehat u}{\partial\xi^t\partial\xi^s}
+\sum_{s=1}^3 b_s(\boldsymbol\xi)\,
\frac{\partial\widehat u}{\partial\xi^s},
\label{eq:laplacepullback}
\end{equation}
where the coefficients are
\begin{equation}
a_{st}(\boldsymbol\xi)=\sum_{r=1}^3 J^{-1}_{sr}(\boldsymbol\xi) J^{-1}_{tr}(\boldsymbol\xi),
\qquad
b_s(\boldsymbol\xi)=\sum_{q,r=1}^3  J^{-1}_{qr}(\boldsymbol\xi)\,
\partial_q(J^{-1}_{sr} (\boldsymbol\xi)  ).
\label{eq:coeffA}
\end{equation}
Rearranging into matrix form yields
\(\boldsymbol A(\boldsymbol\xi)=(a_{st}(\boldsymbol\xi))_{s,t=1}^3=\boldsymbol J^{-1}\boldsymbol J^{-\mathrm T}=(\boldsymbol J^{\mathrm T}\boldsymbol J)^{-1}\).

Thus, \eqref{eq:mappde} is equivalent to
\begin{equation}
-\sum_{s,t=1}^3 a_{st}(\boldsymbol\xi)
\frac{\partial^2\widehat u}{\partial\xi^t\partial\xi^s}
-\sum_{s=1}^3 b_s(\boldsymbol\xi)
\frac{\partial\widehat u}{\partial\xi^s}
=\widehat f(\boldsymbol\xi),
\qquad \boldsymbol\xi\in(0,1)^3,
\label{eq:pullbackPDE}
\end{equation}
and the right-hand term is
\begin{equation}
\widehat f(\boldsymbol\xi)=f\bigl(\boldsymbol G(\boldsymbol\xi)\bigr).
\label{eq:fhatdef}
\end{equation}
At this point, the original problem is equivalently represented as a second-order equation with variable coefficients over the parameter field, where the influence of the geometry mapping is only introduced into the system through the coefficients
\(a_{st}\), \(b_s\) and the right-hand term \(\widehat f\).

The quantity \(b_s\) in \eqref{eq:coeffA} involves parametric derivatives of \(\boldsymbol J^{-1}\). From 
\[
(\partial_q\boldsymbol J)\boldsymbol J^{-1}+\boldsymbol J\,\partial_q(\boldsymbol J^{-1})
=\partial_q(\boldsymbol J\boldsymbol J^{-1})=\partial_q \boldsymbol I=\boldsymbol 0,
\]
we obtain
\[
\partial_q(\boldsymbol J^{-1})=-\boldsymbol J^{-1}(\partial_q\boldsymbol J)\boldsymbol J^{-1}.
\]
Expressed in component form, this equation becomes
\begin{equation}
\partial_q (J^{-1}_{sr}(\boldsymbol\xi))
=-\sum_{k,l=1}^3  J^{-1}_{sl}(\boldsymbol\xi)\,(\partial_q\partial_k G_l)\, J^{-1}_{kr}(\boldsymbol\xi),
\qquad q,s,r=1,2,3.
\label{eq:ttiga:dJ}
\end{equation}

Substituting \eqref{eq:ttiga:dJ} into \eqref{eq:coeffA} and using the contraction
\(\sum_{r=1}^3 J^{-1}_{qr}(\boldsymbol\xi)J^{-1}_{kr}(\boldsymbol\xi)=a_{qk}(\boldsymbol\xi)\), we obtain
\begin{equation}
b_s=-\sum_{q,k,l=1}^3 a_{qk}\,(\partial_q\partial_kG_l)\, J^{-1}_{sl}
=-\bigl(\boldsymbol J^{-1} \boldsymbol h\bigr)_s,
\qquad
h_l=\sum_{q,k=1}^3 a_{qk}\,\partial_q\partial_kG_l .
\label{eq:ttiga:b2}
\end{equation}
where \(\boldsymbol h=(h_l)_{l=1}^3\). Accordingly, \(a_{st}\) and \(b_s\) only require the first- and second-order derivatives of \(\boldsymbol G\), which can be evaluated directly at collocation points.


\subsection{QTT Representation of Discrete Equations in Parameter Domain}

Consistent with the assumptions in Section~3.1, the same number of basis
functions \(n\), polynomial degree \(p\), and knot vector \(\Xi\) are used in
each direction. We set \(n=2^m\),  and the discrete solution is
represented as
\begin{equation}
\widehat u_h(\boldsymbol{\xi})
=
\sum_{i_1,i_2,i_3=1}^n
c_{i_1i_2 i_3}
\prod_{k=1}^3 N_{i_k,p}(\xi^k),
\label{eq:maph}
\end{equation}
where \(c_{i_1 i_2 i_3}\) are the unknown coefficients. The associated
tensor is denoted by \(\boldsymbol C=(c_{i_1i_2 i_3})\) and is
vectorized as \(\operatorname{vec}(\boldsymbol C)\) in the ordering specified
by \eqref{eq:vecC}. Thus,
\(\operatorname{vec}(\boldsymbol C)(\iota(\mathbf i))=c_{i_1 i_2 i_3}\),
where
\(\iota(\mathbf i)=i_3+\sum_{k=1}^{2}(i_k-1)n^{3-k}\).  The collocation points are given by
\begin{equation}
 \boldsymbol{\bar\xi}_{\mathbf i}
=
\bigl(\bar\xi_{i_1},\bar\xi_{i_2},\bar\xi_{i_3}\bigr),
\qquad
\mathbf i=(i_1,i_2,i_3) \in\{1,\dots,n\}^3,
\label{eq:mappoints}
\end{equation}
where the \(\bar\xi_i\) are Greville abscissae.
Since the one-dimensional matrices depend only on the B-splines on the parameter domain and on the collocation points, the one-dimensional matrices \(\boldsymbol D^{(0)},\boldsymbol D^{(1)},\boldsymbol D^{(2)}\) defined in Section~2.1 remain applicable here.

Substituting \eqref{eq:maph} into \eqref{eq:pullbackPDE} and applying the same
analysis as in Section~3.1 yields the discrete matrix associated with the
operator in \eqref{eq:pullbackPDE}. For \(s,t=1,2,3\), define
\begin{equation}
\boldsymbol Q_{st}
:=
\bigotimes_{k=1}^{3}\boldsymbol Q^{(st)}_k,
\label{eq:Pstdef}
\end{equation}
where
\begin{equation}
\boldsymbol Q^{(st)}_k=
\begin{cases}
\boldsymbol D^{(2)}, & s=t\ \text{and}\ k=s,\\
\boldsymbol D^{(1)}, & s\ne t\ \text{and}\ k\in\{s,t\},\\
\boldsymbol D^{(0)}, & \text{otherwise}.
\end{cases}
\label{eq:Akdef}
\end{equation}
We also define
\begin{equation}
\boldsymbol H_s
=
(\boldsymbol D^{(0)})^{\otimes(s-1)}\otimes \boldsymbol D^{(1)}\otimes( \boldsymbol D^{(0)})^{\otimes(3-s)}.
\label{eq:Hsdef}
\end{equation}
The operator on the left-hand side of \eqref{eq:pullbackPDE} is then
discretized as
\begin{equation}
\begin{aligned}
\boldsymbol L
={}
-\sum_{s,t=1}^3
\operatorname{diag}(a_{st})\,\boldsymbol Q_{st}
-\sum_{s=1}^3
\operatorname{diag}(b_s)\,\boldsymbol H_s,
\end{aligned}
\label{eq:LGfull}
\end{equation}
where \(\operatorname{diag}(\cdot)\) denotes the diagonal matrix whose diagonal
entries are the samples at all collocation points, ordered according
to the same multi-index ordering as \(\operatorname{vec}(\boldsymbol C)\),
\begin{equation}
\operatorname{diag}(a_{st})({\iota(\mathbf i),\,\iota(\mathbf i)})
=a_{st}\bigl(\boldsymbol{\bar\xi}_{\mathbf i}\bigr),
\qquad \mathbf i\in\{1,\dots,n\}^3 .
\label{eq:Diagdef}
\end{equation}
The coefficient function \(b_s\) is also arranged in this order. The
right-hand side is discretized as a vector in the same ordering of collocation points,
\begin{equation}
\boldsymbol f ({\iota(\mathbf i))}
=
\widehat f\bigl(\boldsymbol{\bar\xi}_{\mathbf i}\bigr)
=
f\bigl(\boldsymbol G(\boldsymbol{\bar\xi}_{\mathbf i})\bigr).
\label{eq:Fmap}
\end{equation}
When \( \boldsymbol G=\mathrm{\boldsymbol id}\), we have
\(\boldsymbol J=\boldsymbol I\) and hence
\(\boldsymbol J^{-1}=\boldsymbol I\). Equation~\eqref{eq:coeffA} then gives
\(a_{st}=\delta_{st}\) and \(b_s=0\). In this case, only the second-order terms
with \(s=t\) remain in \eqref{eq:LGfull}, and the operator reduces to the
standard discretized Laplace operator given by \eqref{eq:Lfull}.

Section~3.2.2 provides the QTT-matrix representations of
\(\boldsymbol D^{(0)}\), \(\boldsymbol D^{(1)}\), and
\(\boldsymbol D^{(2)}\). Consequently, the QTT-matrix representations of \(\boldsymbol Q_{st}\) and \(\boldsymbol H_s\) in \eqref{eq:Pstdef} and \eqref{eq:Hsdef} are directly obtained from the Kronecker product in the TT format, as in \eqref{eq:tt_Kronecker}.
To construct the QTT representations of multidimensional vectors containing the values of \(a_{st}\) and \(b_s\)
at all collocation points, we use the TT-cross interpolation. The elements of these multi-dimensional vectors are not all assembled, but are sampled by evaluating functions on carefully selected multi-index entries. This process is guided by a recursive rank-adaptive algorithm based on the principle of maximum volume.  The corresponding diagonal QTT matrices \(\operatorname{diag}(a_{st})\), and
\(\operatorname{diag}(b_s)\) are then
constructed from these QTT vectors as in \eqref{eq:diag_QTT}. The QTT vector for the right-hand
side \(\boldsymbol f\) is obtained in the same manner by the TT-cross interpolation.
The QTT matrices for \(\boldsymbol Q_{st}\),
\(\boldsymbol H_s\), \(\operatorname{diag}(a_{st})\), and
\(\operatorname{diag}(b_s)\) are assembled termwise according to
\eqref{eq:LGfull}, using the TT contraction in
\eqref{eq:tt_local_matvec} and the TT summation in \eqref{eq:tt_sum},
followed by TT rounding. Then the QTT representation of
\(\boldsymbol L\) is obtained.

The  handling  of the boundary conditions is completely analogous to that in Section 3.3 and is therefore omitted for brevity.

\section{Numerical Experiments}

This section examines the QTT-based collocation solver from three perspectives. First of all,  the one-dimensional experiment verifies whether the QTT representation of the
banded collocation matrices proposed in Section 3.2 maintains stable ranks under mesh refinement. Secondly, in order to verify the QTT-based solver under
identity geometry mapping as described in Section 3, we solve the convection-diffusion-reaction problem as the spatial dimension increases. Finally, by performing mesh refinement on annular-region and shell examples, we evaluate the performance of the QTT-based solver for solving 3D Poisson problems under general geometry mappings proposed in Section 4.

In this experiment, tensor operations are performed with \texttt{ttpy} \cite{oseledets2013ttp}.
TT-cross interpolation of the variable-coefficient fields and multivariable functions is implemented by \texttt{tt.multifuncrs}. The final algebraic systems for the QTT-based method is 
\begin{equation}
    \mathbf L\mathbf c=\mathbf f,
    \label{eq:qtt-linear-system}
\end{equation}
where \(\mathbf L\) is the discrete operator, \(\mathbf c\) is the unknown B-spline
coefficients, and \(\mathbf f\) is the final right-hand side. The system is solved for the
unknown QTT vector \(\mathbf c\) by \texttt{tt.amen\_solve}, which updates one
coefficient QTT core at a time to reduce the residual. As a comparative experiment, when solving using the classical full-grid method, the algebraic system is expressed as
\(\boldsymbol{L}\boldsymbol{c}=\boldsymbol{f}\).  In the experiment, we solve this full-grid system by the unpreconditioned GMRES method.

\subsection{QTT Structure of the 1D Collocation Matrices}
\label{sec:one-d}

The purpose of this experiment is to solve a partial differential equation,
but to verify whether the one-dimensional collocation matrices can be accurately
represented by the QTT representation. For the numbers of basis functions
\(n=16,32,\ldots,4096\) and the B-spline degrees \(p=2,3,\ldots,7\), the QTT representations of matrices \(\boldsymbol D^{(0)}, \boldsymbol D^{(1)}, \boldsymbol D^{(2)}\) are constructed and compared with the corresponding dense $n \times n$ matrices. The tolerance of TT rounding used in the banded-vector construction is \(10^{-13}\). For a dense matrix \(\boldsymbol A\), the relative Frobenius error between \(\boldsymbol A\) and its QTT representation is defined as
\begin{equation}
    \frac{\|\boldsymbol A_{\mathrm{QTT}}-\boldsymbol A\|_F}{\|\boldsymbol A\|_F},
    \label{eq:frobenius-error}
\end{equation}
where \(\boldsymbol A_{\mathrm{QTT}}\) denotes the matrix obtained by expanding the QTT representation of \(\boldsymbol A\) into its dense matrix form.

\begin{figure}[H]
\centering
\includegraphics[width=\textwidth]{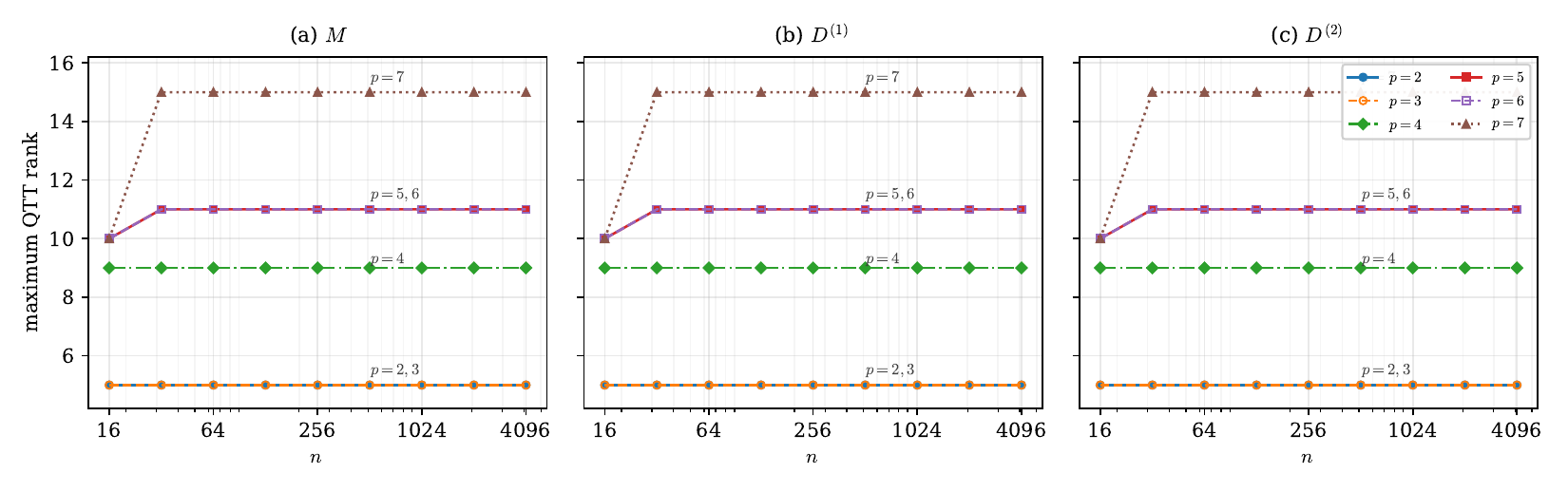}
\caption{Maximum QTT ranks of the one-dimensional collocation matrices versus the number of basis functions
\(n\). The three subfigures correspond to \(\boldsymbol D^{(0)}\), \(\boldsymbol D^{(1)}\), and \(\boldsymbol D^{(2)}\), and
the curves correspond to different B-spline degrees \(p\).}
\label{fig:one-d-rank}
\end{figure}

Figure~\ref{fig:one-d-rank} presents the maximum QTT ranks of \(\boldsymbol D^{(0)}\),
\(\boldsymbol D^{(1)}\), and \(\boldsymbol D^{(2)}\).
For different $p$-values, the ranks show a clear saturation as $n$ increases. For
\(p=2,3\), the maximum rank of all three matrices is 5 over the whole mesh
range. For \(p=4\), the maximum rank is 9. For \(p=5,6\), the maximum rank is
11 except at \(n=16\). For \(p=7\), the rank is 10 at \(n=16\) and becomes 15
for \(n\ge32\). After the first mesh refinement, each curve becomes almost
horizontal. These results show that the maximum
QTT rank is mainly determined by the B-spline degree \(p\), which is consistent with our discussion in subsection 3.2.2.

\begin{figure}[H]
\centering
\includegraphics[width=\textwidth]{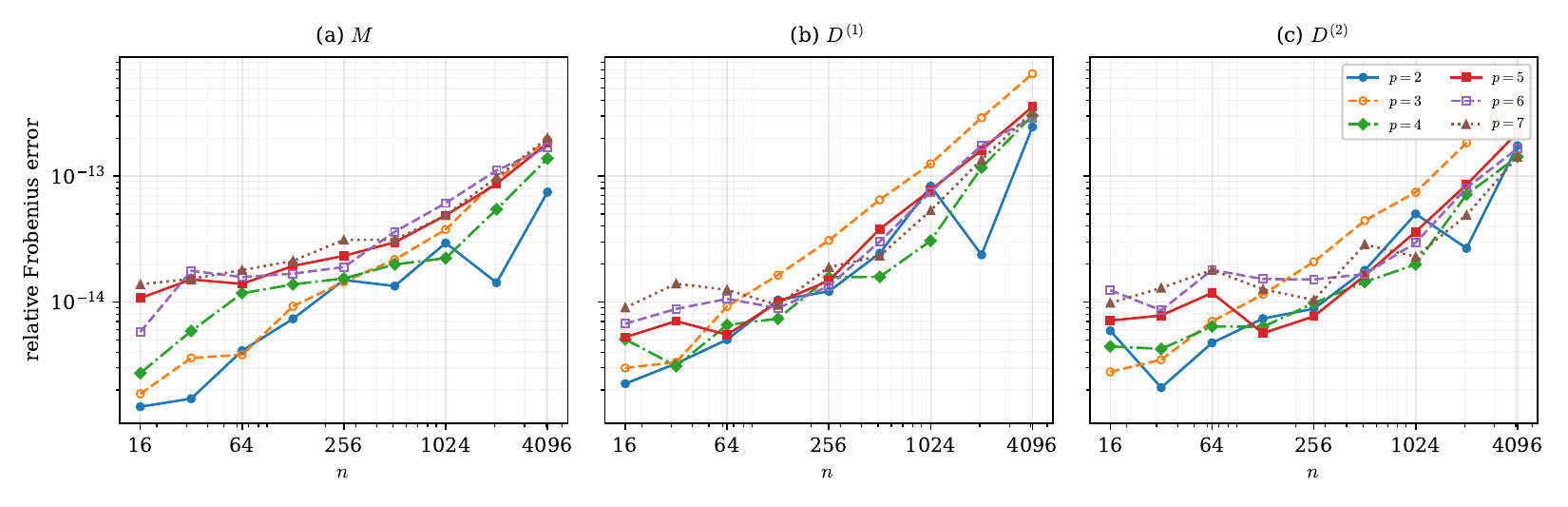}
\caption{The relative Frobenius errors between the matrices \(\boldsymbol D^{(0)}\), \(\boldsymbol D^{(1)}\),
\(\boldsymbol D^{(2)}\) and their QTT representations versus the number of basis functions
\(n\). The errors are computed from
\eqref{eq:frobenius-error} and
the curves correspond to different B-spline degrees \(p\).}
\label{fig:one-d-error}
\end{figure}

Figure~\ref{fig:one-d-error} presents the relative Frobenius errors of 
\(\boldsymbol D^{(0)}\),
\(\boldsymbol D^{(1)}\), and \(\boldsymbol D^{(2)}\) defined in \eqref{eq:frobenius-error}. The errors for all \(p\) remain close
to the  construction tolerance \(10^{-13}\). As \(n\) increases, the error will slightly increase,
which is mainly due to the TT rounding. The QTT representation
therefore reproduces \(\boldsymbol D^{(0)},\boldsymbol D^{(1)},\boldsymbol D^{(2)}\) accurately, and its rank and storage
do not grow exponentially with the number of basis functions.

\subsection{Convection-Diffusion-Reaction Problem  on regular Cartesian domains}
\label{sec:steady}

This example examines the performance of the QTT-based solver as the spatial dimension $d$ increases. The domain is the unit hypercube \(\Omega=(0,1)^d\), and
the problem is
\begin{equation}
    -\Delta u(\boldsymbol x)+\boldsymbol b(\boldsymbol x)\cdot\nabla u(\boldsymbol x)
    +c(\boldsymbol x)u(\boldsymbol x)=f(\boldsymbol x),
    \qquad
    u(\boldsymbol x)|_{\partial\Omega}=g(\boldsymbol x).
\end{equation}
The analytical solution is
\begin{equation}
    u_{\mathrm{ex}}(\boldsymbol x)   =\sin\left(2\pi\sum_{s=1}^{d}x_s\right),\quad\boldsymbol x=(x_1,\ldots,x_d).\label{eq:CDRexact}
\end{equation}
The variable
coefficients $\boldsymbol b(\boldsymbol x)$ and $c(\boldsymbol x)$ are
\begin{equation}  
    b_s(x_s)=
    \begin{cases}
    \sin(2\pi x_s),&s\text{ is odd},\\
    \cos(2\pi x_s),&s\text{ is even},
    \end{cases}\qquad
    c(\boldsymbol x)=\cos\left(2\pi\sum_{s=1}^{d}x_s\right).
\end{equation}
The right-hand side force function $f(\boldsymbol{x})$ and the boundary conditions $g(\boldsymbol{x})$ are determined in
accordance with the exact solution \eqref{eq:CDRexact}.

The discretization uses \(n=128\) and \(p=2\) in every direction, with
\(d=1,\ldots,8\). The diffusion term is assembled according to \eqref{eq:Lblock}- \eqref{eq:block_qtt_col}. The coefficients \(b_s,c\) and the function
\(f\) are constructed by TT-cross with tolerance \(10^{-10}\).  The classical full-grid system is computed up to \(d=3\) by the 
unpreconditioned GMRES method. The larger dimension exceeds the computational capability of our equipment.
The relative \(L^2\) error against the analytical solution is
\begin{equation}
    \sqrt{
    \frac{\sum_{\mathbf i}(\hat u_{h}(\boldsymbol{x}_\mathbf i)-u_{\mathrm{ex}}(\boldsymbol{x}_\mathbf i))^2}
    {\sum_{\mathbf i}u_{\mathrm{ex}}(\boldsymbol{x}_\mathbf i)^2}
    },
    \label{eq:common-error}
\end{equation}
where \(\hat u_h\) denotes the numerical solution.
In this subsection, the multi-index \(\mathbf i\) runs over the complete \(9^d\) tensor grid
formed by nine uniform points in each direction. 

The reported CPU time starts from the assembly of one-dimensional QTT representations, and ends with the linear solution and the final TT rounding. The evaluation of the analytical error is excluded. To quantify the compression achieved, we introduce
the compression ratio, defined as
\begin{equation}
    \frac{N_{\mathbf L}}{\stor(\mathbf L)},
    \qquad
    \frac{n^d}{\stor(\mathbf c)},
    \qquad
    \frac{n^d}{\stor(\mathbf f)},
    \label{eq:common-compression}
\end{equation}
where \(\stor(\mathbf L),\stor(\mathbf c)\), and \(\stor(\mathbf f)\) denote
the actual QTT storage of the respective objects. The reference size of $\mathbf L$ is
\(
    N_{\mathbf L}=(5n-6)^d,
\)
the number of nonzero elements in a \(d\)-dimensional operator formed as the tensor product of  one-dimensional five-band collocation matrices.

\begin{figure}[H]
\centering
\includegraphics[width=\textwidth]{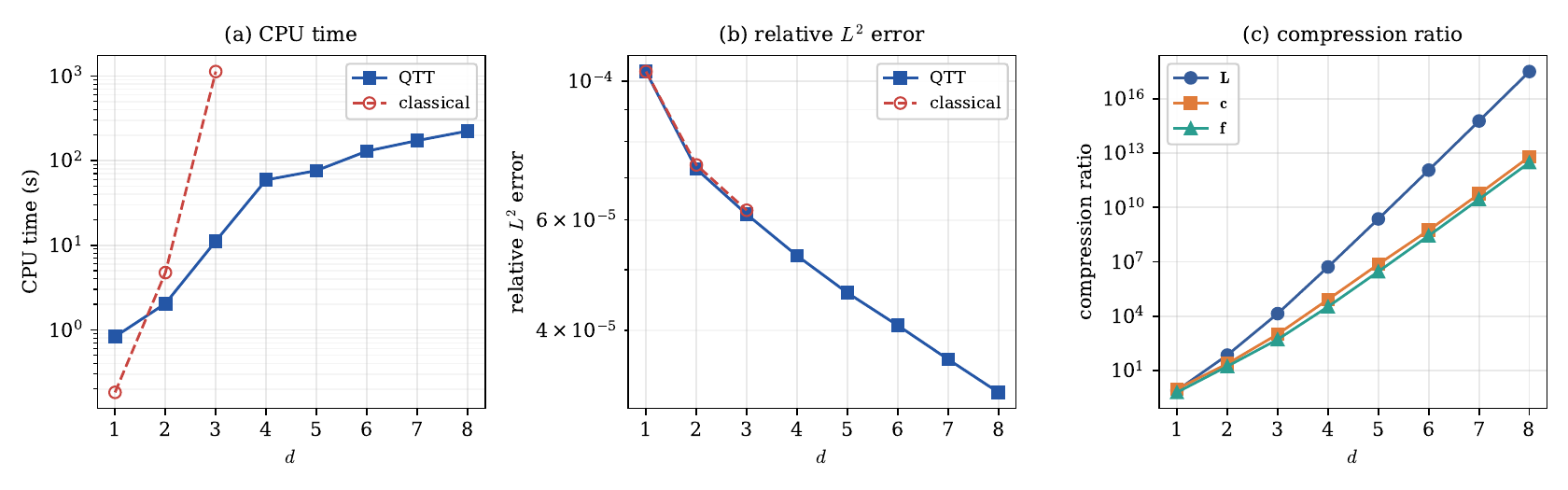}
\caption{Performance comparison of the classical full-grid method and the QTT-based method for the convection-diffusion-reaction problem at different spatial dimensions $d$.  (a) shows CPU time. (b) shows the relative \(L^2\) error defined in \eqref{eq:common-error}.   (c) shows the compression
ratios of \(\mathbf L,\,\mathbf c\) and \(\mathbf f\) defined in \eqref{eq:common-compression}.}
\label{fig:steady-overview}
\end{figure}

Figure~\ref{fig:steady-overview} presents the CPU time, the relative \(L^2\) error, and the compression ratio at different spatial dimensions \(d\) for the classical full-grid method and the QTT-based method. In subfigure~\ref{fig:steady-overview}(a), the two methods show different growth patterns of the CPU time. the CPU time of the QTT-based method slowly increases from 0.83 s at \(d=1\) to 11.19 s at \(d=3\), and reaches
223.39 s at \(d=8\). The CPU time of the classical full-grid method rapidly increases from 0.18 s at \(d=1\) to 1128.24 s at \(d=3\). 
For $d=1$ and $n=128$, the matrix-free classical method has a simple process flow and very little overhead, and is slightly faster than the method based on the QTT-based method. The QTT-based method is faster than the classical full-grid method for \(d>1\), and about 101 times faster at \(d=3\).

Figure~\ref{fig:steady-overview}(b) shows the errors of the two methods relative to the analytical solution defined in \eqref{eq:common-error}. 
For \(d=1,2,3\), the error curves of the classical full-grid method and the QTT-based method nearly coincide, indicating that the QTT-based method  achieves its speedup without sacrifice accuracy. With a further increase in dimension, the error of the QTT-based method varies only slightly.

Figure~\ref{fig:steady-overview}(c) the compression
ratios of \(\mathbf L,\,\mathbf c\) and \(\mathbf f\) defined in \eqref{eq:common-compression}
at different spatial dimensions $d$. As $d$ increases, the compression ratios of $\mathbf L$, $\mathbf c$, and $\mathbf f$ also rapidly increase. Our QTT-based strategy significantly reduces memory usage. This shows that,
as the dimension increases, QTT maintains almost the same analytical accuracy
as the classical full-grid method while avoiding the high storage cost of the complete
high-dimensional tensor and reducing the CPU time.

\subsection{Poisson Problems on
3D nontrivial geometries}

The previous example mainly examine the performance of our QTT-based algorithm on Cartesian grids as the spatial dimension varies. To assess the method in non-Cartesian domains, we consider Poisson problems \eqref{eq:poisson} on a thick quarter of annulus domain, for which \(\boldsymbol J(\boldsymbol{\xi})\) in \eqref{eq:jacdef} is diagonal, and a spherical shell, for which \(\boldsymbol J(\boldsymbol{\xi})\) has nine nonzero entries. Both examples use \(p=2\) and
\(n=16,32,64,128,256,512\) in the QTT-based method. The classical full-grid method is computed up to \(n=128\), since the computations with \(n>128\) exceeded the computing capacity available for this study.

Unless stated otherwise, the TT-cross parameters, error definitions, timing
convention, and compression definitions are the same as in
Sections~\ref{sec:one-d} and~\ref{sec:steady}. The error against the analytical
solution is computed
from \eqref{eq:common-error} on the complete \(n^3\) tensor grid formed by
\(n\) uniform points in each parameter direction.

\subsubsection{Thick Quarter of Annulus Domain}
\label{sec:ring}

The geometry mapping is
\begin{equation*}
    \boldsymbol G(\xi,\eta,\zeta)
    =
    \bigl(\rho(\xi)\cos\theta(\eta),
    \rho(\xi)\sin\theta(\eta),\zeta\bigr),
\end{equation*}
with
\begin{equation*}
    \rho(\xi)=r_{\rm in}+(r_{\rm out}-r_{\rm in})\xi,
    \qquad
    \theta(\eta)=\frac{\pi}{2}\eta,
\end{equation*}
where \(r_{\rm in}=0.5\) and \(r_{\rm out}=1\). Dirichlet conditions
\(u(r_{\rm in})=1\) and \(u(r_{\rm out})=2\) are imposed on the two cylindrical
surfaces, and the other four surfaces have homogeneous Neumann conditions. This is a homogeneous Laplace problem with $f= 0$. The
exact solution depends only on the radius,
\begin{equation}
    u_{\mathrm{ex}}(r)=
    \frac{
    u_{\rm in}\log(r_{\rm out}/r)
    +u_{\rm out}\log(r/r_{\rm in})
    }{
    \log(r_{\rm out}/r_{\rm in})
    }.
    \label{eq:ring-exact}
\end{equation}

 Only the coefficient functions \(a_{11},a_{22},a_{33},b_1,b_2,b_3\) are nonzero, and their QTT representations are obtained using the TT-cross settings of Section~\ref{sec:steady}.
Two assembly errors are defined,
\begin{equation}
    \frac{ 
    \| \operatorname{vec}( \mathbf L\mathbf v)-Lv\|
    }{
    \|Lv\|
    },
    \qquad
    \frac{\|\operatorname{vec}(\mathbf f)-f\|}{\|f\|},
    \label{eq:assembly-error}
\end{equation}
where \(\mathbf v\) is a random TT vector, \(v\) is its full-grid vector.

\begin{figure}[H]
\centering
\includegraphics[width=\textwidth]{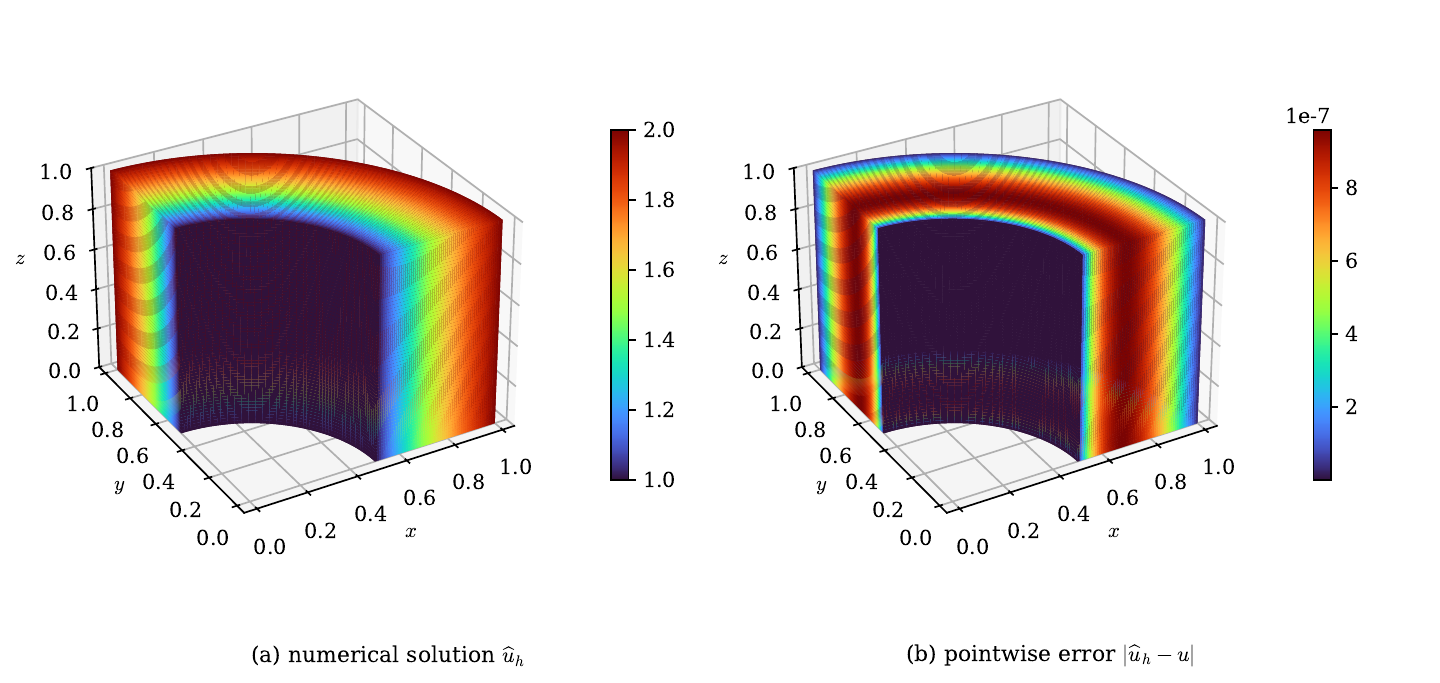}
\caption{Thick quarter of annulus domain with \(n=128\) and \(p=2\). (a)
shows the QTT numerical solution and (b) shows the pointwise
absolute error.}
\label{fig:ring-solution}
\end{figure}

\begin{figure}[H]
\centering
\includegraphics[width=0.92\textwidth]{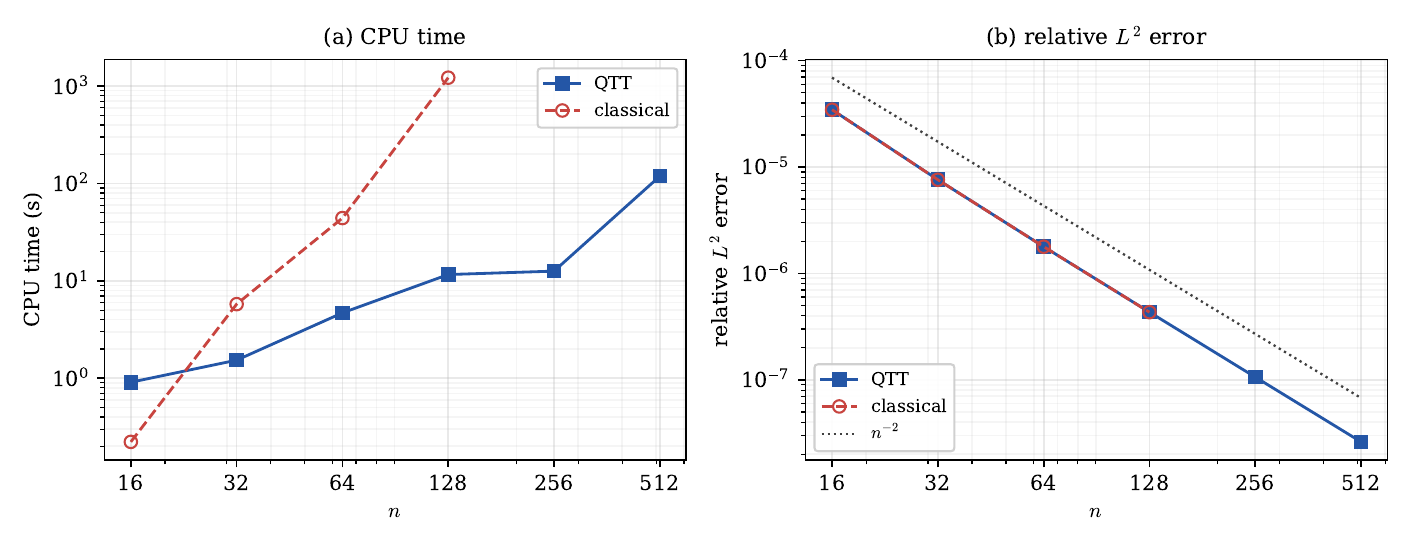}
\caption{Performance comparison of the classical full-grid method and the QTT-based method in the thick quarter of annulus domain. (a)
shows CPU time and (b) shows the relative \(L^2\) error defined in \eqref{eq:common-error} versus the number of basis functions \(n\) per direction.  The dotted line gives the
\(n^{-2}\) reference slope.}
\label{fig:ring-performance}
\end{figure}

Figure~\ref{fig:ring-solution} shows the QTT-baesd  numerical solution and the corresponding pointwise absolute error at \(n=128\). Figure~\ref{fig:ring-performance} compares the performance of the QTT-based method and the classical method in terms of CPU time and the error against the analytical
solution .  In subfigure~\ref{fig:ring-performance}(a),  for $n = 16$, the QTT-based method performs more slowly than the classical method. The QTT-based method becomes faster than the classical method starting from \(n=32\), and at \(n=128\) it is approximately 105 times faster.
In subfigure~\ref{fig:ring-performance}(b), the QTT and classical error curves nearly coincide over the
common range. The QTT relative \(L^2\) error decreases from
\(3.457\times10^{-5}\) at \(n=16\) to \(2.613\times10^{-8}\) at \(n=512\).
The fitted observed order is about 2.07, parallel to the \(n^{-2}\) reference
line. The \(n^{-2}\) reference line represents the expected second-order convergence of the quadratic B-spline collocation discretization \cite{math11020469}, since \(h=O(n^{-1})\) and the discretization error is \(O(h^{2})\). 

\begin{figure}[H]
\centering
\includegraphics[width=0.92\textwidth]{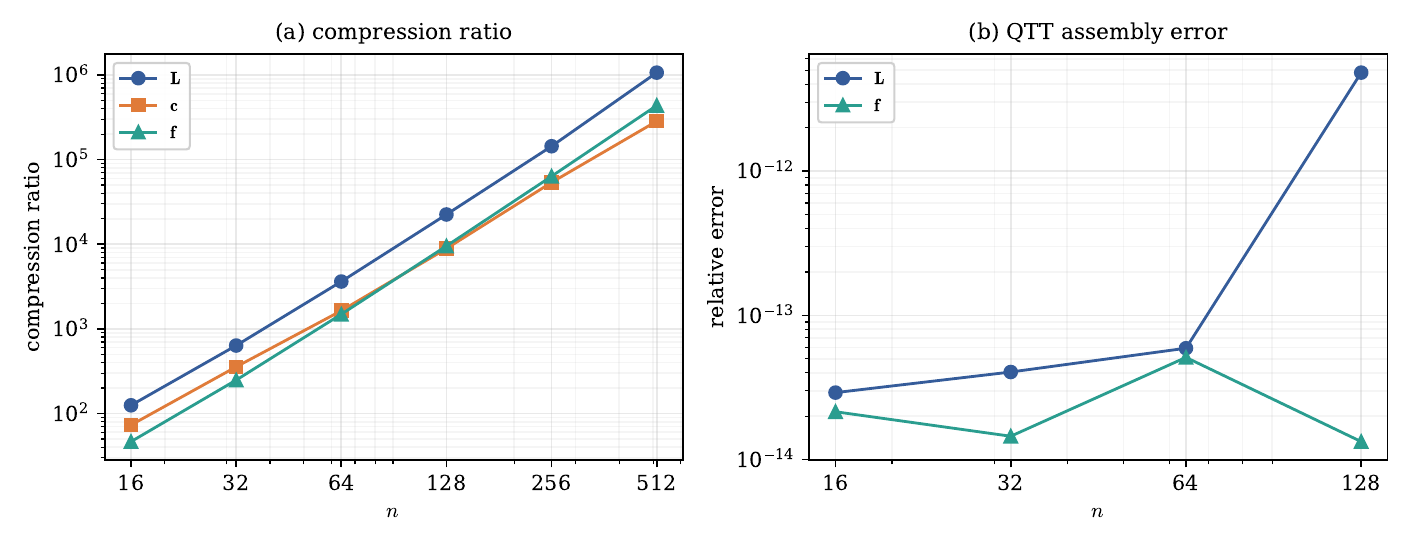}
\caption{Storage and assembly errors for the thick quarter of annulus domain at different numbers of basis functions \(n\) per direction. (a)
shows the compression ratios of \(\mathbf L,\mathbf c,\mathbf f\) in
\(\mathbf L\mathbf c=\mathbf f\), defined in
\eqref{eq:common-compression}.  (b) shows the QTT operator error $\mathbf L$ and right-hand side error $\mathbf f$, defined in
\eqref{eq:assembly-error}.}
\label{fig:ring-compression}
\end{figure}

Figure~\ref{fig:ring-compression} shows the storage and assembly errors defined in
\eqref{eq:assembly-error}. In
subfigure~\ref{fig:ring-compression}(a), the compression ratios of \(\mathbf{L}\), \(\mathbf{c}\), and \(\mathbf{f}\) increase steadily with \(n\). When $n=512$, the compression ratios of \(\mathbf L,\mathbf c,\mathbf f\) are between $10^{5}$ and $10^6$. In subfigure~\ref{fig:ring-compression}(b), for \(n\le128\), the
operator and right-hand side errors defined in
\eqref{eq:assembly-error} are much smaller than the
error in subfigure~\ref{fig:ring-compression}(b) throughout the same range. Therefore, the global error
 is still controlled by
spatial discretization.

For this annulus example, as the number of basis functions per direction, the QTT method keeps the solution accuracy almost unchanged
while greatly reducing the CPU time and increasing the compression ratio. The
benefit becomes more substantial as \(n\) increases.

\subsubsection{Spherical Shell}
\label{sec:shell}

The shell uses an ordinary B-spline approximation with a \(2\times3\times3\)
control net and unit weights as in \cite{mantzaflaris2017}. Its radial thickness is
\(r_{\rm out}-r_{\rm in}=0.2\), with \(r_{\rm in}=0.8\) and
\(r_{\rm out}=1.0\). Unlike the ring, all six metric components \(a_{st}\) and
all three coefficients \(b_s\) are nonzero.
The manufactured solution and source are
\begin{equation}
    u_{\rm ex}(x,y,z)=\sin(\pi x)\sin(\pi y)\sin(\pi z),
    \qquad
    f=3\pi^2u_{\rm ex}.
\end{equation}
Dirichlet boundary data consistent with the manufactured solution are imposed
on all six faces. 

\begin{figure}[H]
\centering
\includegraphics[width=\textwidth]{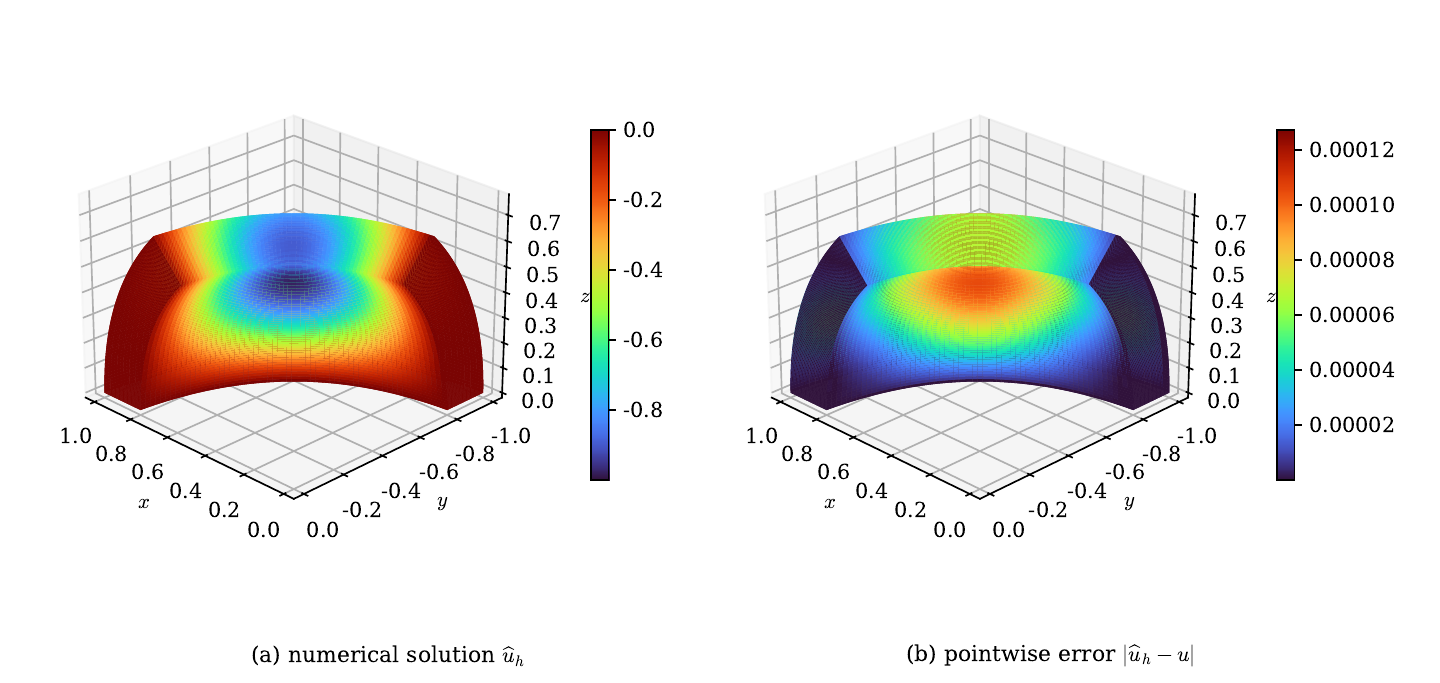}
\caption{Quarter-spherical shell with \(n=128\) and \(p=2\). (a)
shows the QTT numerical solution and (b) shows the pointwise
absolute error.}
\label{fig:shell-solution}
\end{figure}

\begin{figure}[H]
\centering
\includegraphics[width=0.92\textwidth]{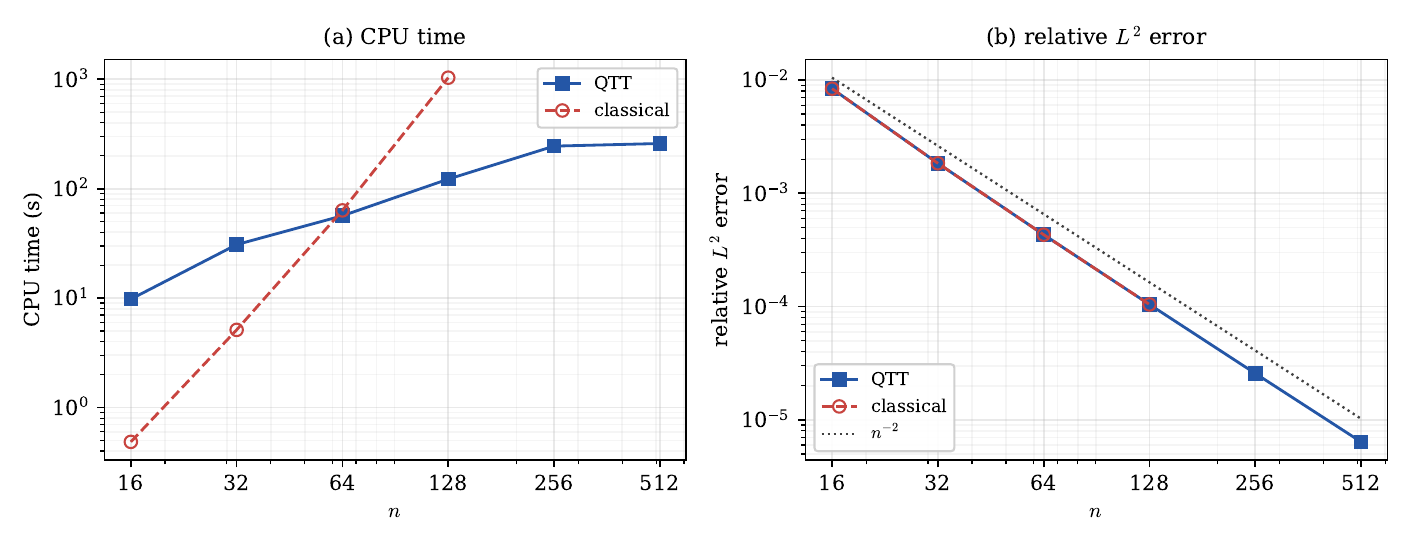}
\caption{Performance comparison of the classical full-grid method and the QTT-based method in the sphere shell. (a) shows CPU time
and (b) shows the relative \(L^2\) error defined in \eqref{eq:common-error} versus the number of basis functions
\(n\) per direction. The dotted line gives the \(n^{-2}\)
reference slope.}
\label{fig:shell-performance}
\end{figure}

Figure~\ref{fig:shell-solution} shows the QTT-baesd numerical solution and the corresponding pointwise abso-
lute error  at
\(n=128\). Figure~\ref{fig:shell-performance} compares the performance of the QTT-based method and the
classical method in terms of CPU time and the error against the analytical
solution . In subfigure~\ref{fig:shell-performance}(a),  the classical method increases sharply from
0.49 s at \(n=16\) to 1033.37 s at \(n=128\), while the
QTT-based method slowly increases from
9.79 s at \(n=16\) to 258.67 s at \(n=512\). The two curves cross near
\(n=64\). At \(n=128\), the
QTT-based method is about 8.5 times faster.
In subfigure~\ref{fig:shell-performance}(b), the relative \(L^2\) error of the
QTT-based method decreases from
\(8.388\times10^{-3}\) at \(n=16\) to \(6.394\times10^{-6}\) at \(n=512\).
The classical curve nearly coincides with the QTT curve for \(n\le128\). The
fitted observed order is about 2.07, parallel to the \(n^{-2}\) reference line.

\begin{figure}[H]
\centering
\includegraphics[width=0.92\textwidth]{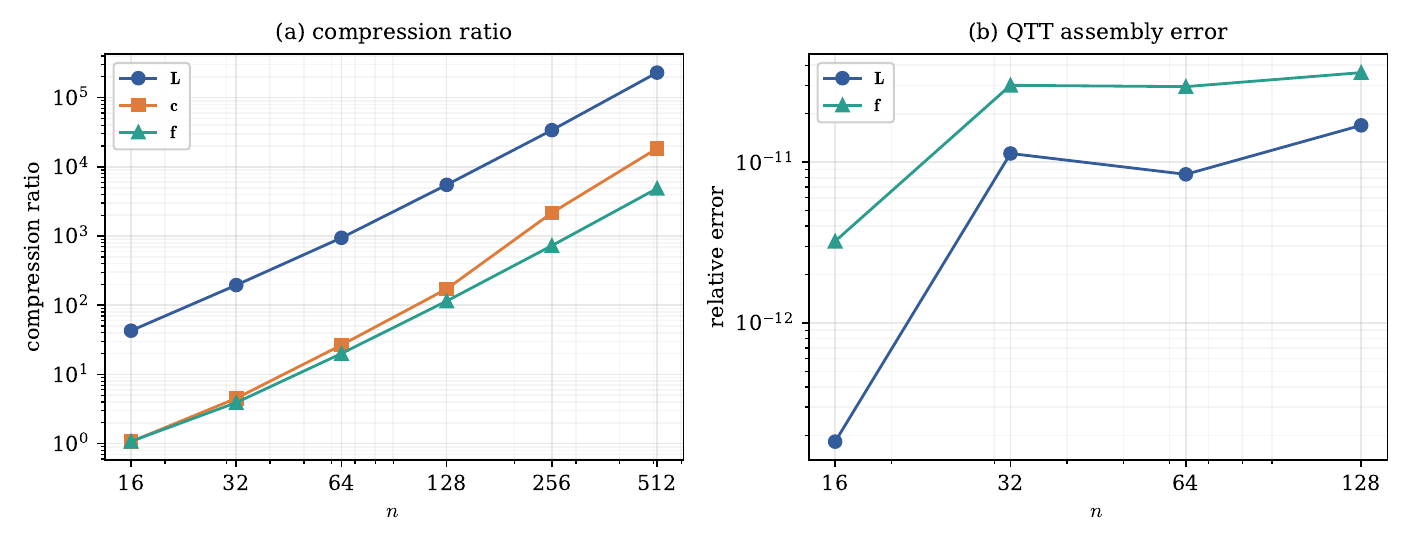}
\caption{Storage and assembly errors for the quarter-spherical shell at different numbers of basis functions   \(n\) per direction.  (a)
shows the compression ratios of \(\mathbf L,\mathbf c,\mathbf f\) in
\(\mathbf L\mathbf c=\mathbf f\), defined in
\eqref{eq:common-compression}.  (b) shows the QTT operator-probe error
of $\mathbf L$ and right-hand side error $\mathbf f$, defined in
\eqref{eq:assembly-error}.}
\label{fig:shell-compression}
\end{figure}

Figure~\ref{fig:shell-compression} shows the storage and assembly errors. In
subfigure~\ref{fig:shell-compression}(a), all three compression ratios exhibit a steady increasing trend with \(n\). At \(n=512\), the compression ratios of \(\mathbf L\), \(\mathbf c\), and \(\mathbf f\) reach \(2.30\times10^{5}\), \(1.85\times10^{4}\), and \(4.87\times10^{3}\), respectively. In subfigure~\ref{fig:shell-compression}(b), for
\(n\le128\), both the operator and right-hand side errors defined in
\eqref{eq:assembly-error} remain below
\(10^{-10}\).

For this shell example, the QTT-based method preserves the accuracy of the classical full-grid
solver while avoiding the rapid growth of CPU time and storage. Moreover, this advantage becomes more significant as \(n\) increases.

\section{Conclusion}

We have presented a QTT-based isogeometric collocation framework for
solving partial differential equations. The solution space is spanned by tensor-product B-spline bases constructed from open uniform knot vectors, and the Greville abscissae is used as collocation points. It is rigorously shown that each band of these one-dimensional collocation matrices is constant in the interior, except at a few positions near the two ends. This property is used to construct the QTT
representation of the one-dimensional collocation matrices. The required
number of sampled entries and the resulting maximum QTT ranks are at most of
quadratic order of the B-spline degree.

For regular Cartesian domains, the QTT ranks of the high-dimensional
discrete Laplacian are determined by the QTT ranks  of the
one-dimensional collocation matrices and are independent of the spatial
dimension. For nontrivial geometry mappings, we consider the three-dimensional Poisson equation and reformulate it equivalently on the parameter domain.  Geometry-induced effects enter through several variable coefficients,
whose vectors at the collocation points are represented in the QTT format by
TT-cross. The resulting
parameter derivatives allow the same one-dimensional QTT construction to be
used.

Numerical experiments confirm that the maximum QTT rank is governed primarily by the B-spline degree and is essentially independent of the number of basis functions. We then solve the convection–diffusion–reaction problems in different spatial dimensions, thereby validating the capability of our QTT-based isogeometric collocation solver for high-dimensional PDEs. Finally, the solver is tested on Poisson problems over nontrivial
three-dimensional geometries, demonstrating its capability for partial
differential equations on non-Cartesian domains.

The present study is restricted to geometries modeled as a single patch.
Extending the proposed QTT-based isogeometric collocation solver to multi-patch
domains is an important direction for future work.

\section*{Acknowledgements}\label{Acknowledgements}

\end{document}